\documentclass{article}
 \usepackage{amsmath}

\newtheorem{theorem}{Theorem}
\newtheorem{corollary}{Corollary}
\newtheorem{lemma}{Lemma}
\newtheorem{proposition}{Proposition}

\newtheorem{definition}{Definition}
\newtheorem{remark}{Remark}

\def\grad{\nabla}
\def\R{{\rm I}\!{\rm  R}}
\def\la1{\lambda_1}
\def\ph1{\varphi_1}
\newcommand{\N}{{\bf N}}

\title{Eigenvalue, maximum principle and regularity for fully non linear homogeneous operators.}    
\author{I. Birindelli, F. Demengel}
\date{}

\begin{document}
\maketitle

\medskip
{\footnotesize
 \centerline{Dipartimento di Matematica}
  \centerline{Universit\`a di Roma "La Sapienza"}
   \centerline{Piazzale Aldo Moro 5}
   \centerline{00185 Roma, Italia}
} 

\medskip
{\footnotesize
 \centerline{Laboratoire  d'Analyse, G\'eom\'etrie et Modelisation}
  \centerline{Universit\'e de Cergy-Pontoise}
   \centerline{Site de Saint-Martin, 2 Avenue Adolphe Chauvin}
   \centerline{95302 Cergy-Pontoise, France}
} 

\medskip

 \centerline{(Communicated by Aim Sciences)}
 \medskip

\begin{abstract}
The main scope of this article is to define the concept of principal eigenvalue for fully non linear second order operators in bounded domains that are elliptic, homogenous with lower order terms.
In particular we prove maximum and comparison principle, H\"older and Lipschitz regularity. This leads to the existence of a first eigenvalue and eigenfunction and to the existence of solutions of Dirichlet problems within this class of operators.
\end{abstract}
\section{Introduction}
In \cite{BD2}, inspired by the acclaimed work of Berestycki, Nirenberg and Varadhan \cite{BNV}, we extended the definition of principal eigenvalue to Dirichlet problems for fully-non linear second order elliptic operators.

 Precisely,
given a bounded domain
 $\Omega$,  given  $\alpha>-1$  we defined the "principal eigenvalue" for
$F(\grad u,D^2u)$ satisfying:

\begin{itemize}
 \item[(H1)] $F(tp,\mu X)=|t|^{\alpha}\mu F(p,X)$, $\forall t\in
\R^\star$,
$\mu\in\R^+$ 

\item[(H2)] $a |p|^\alpha {\rm tr} N\leq F(p,M+N)-F(p,M)\leq A |p|^\alpha{\rm tr} N$ 
for $0<a\leq A$, $\alpha>-1$ and $N\geq 0$.
\end{itemize}

\noindent Indeed we showed that

$$\bar\lambda=\sup\{ \lambda\in \R,\ \ \exists \  \phi >0 
 \ {\rm in}\
\overline \Omega,\ F(\grad \phi,D^2\phi)+\lambda \phi^{\alpha+1}\leq 0 \mbox{ in
the viscosity sense}
\ \}$$
is well defined and it satisfies the following properties:

\noindent {\bf (I)}{\it There exists $\phi$ a continuous  
positive viscosity solution of 
$$\left\{\begin{array}{lc}
F(\grad \phi,D^2\phi) +\bar\lambda\phi^{\alpha+1}=0 & {\rm in}\ \Omega\\
\phi=0 &  {\rm on}\ \partial \Omega.
\end{array}
\right.
$$ }
{\bf (II)}{\it Furthermore, suppose that $\lambda<\bar\lambda$. If $f\leq 0$ is
bounded and continuous  in $\Omega$, then there exists  $u$ non-negative,  viscosity solution of
\begin{equation}\label{eqif}\left\{\begin{array}{lc}
F(\grad u,D^2u) + \lambda u^{\alpha+1}=f & {\rm in}\ \Omega\\
u=0 &  {\rm on}\ \partial \Omega.
\end{array}
\right.
\end{equation}
If moreover $f$ is negative in $\Omega$, the solution is
unique.          
}

Hence $\bar\lambda$ was denoted {\it principal eigenvalue} of $-F$ in $\Omega$.

In the case $\alpha=0$, and for $F$ a linear uniformly elliptic second order operator, these results are included in  \cite{BNV};  when $F$ is  one of the Pucci  operators the problem has been treated by Quaas \cite{Q} and  Busca, Esteban and Quaas \cite{BEQ}. Their papers give a more complete description of the spectrum and also treat bifurcation problems. In \cite{BD2} and in this note the situation is complicated by the fact that there are no known results about the regularity of the solution, or the existence of the solution, even without the zero order term.

Clearly the operator $F$ can be seen as a non-variational extension of the $p$-Laplacian: 
$\Delta_p=\rm{div}(|\grad .|^{p-2}\grad .)$ with $\alpha=p-2$.

\bigskip
 The scope of this article
is to complete the results of \cite{BD2}; indeed we consider operators that depend explicitly on $x$,  we include lower order terms, moreover we define  $\bar\lambda$ in a more suitable way i.e. without requiring that  super-solutions are positive up to the boundary. Precisely we shall study existence of solutions, eigenvalue problems and regularity of the solutions for operators of the following type:
$$G(x,u, \nabla u, D^2 u):= F(x,\nabla u, D^2u)+ b(x).\nabla u |\nabla
u|^\alpha + c(x)|u|^\alpha u
$$
where $F$ satisfies assumptions as in \cite{BD1} i.e. the above assumption (H1) and (H2), plus some continuity with respect to the $x$ variable. See e.g. \cite{I} for similar conditions. 
Because of the new setting the proofs differ in nature from \cite{BD2} .

The hypothesis  on 
 $b$  and $c$ are quite standard and they  will be described in the next section.

 As mention above we define   $\bar\lambda$ in a more "correct"  way i.e. : 
$$\bar\lambda:=\sup\{ \lambda\in \R,\ \ \exists \  \phi>0\  
\ {\rm in}\ \Omega,\ G(x,\phi, \grad \phi,D^2\phi)+\lambda
\phi^{\alpha+1}\leq 0
\mbox{ in the viscosity sense}
\ \}.$$

\bigskip
The main aim of this paper is to prove the two following  existence results:

\noindent{\it Suppose that $f\leq 0$, bounded and continuous, that $\lambda < \bar\lambda$, then
there exists a non-negative solution of 
$$\left\{\begin{array}{cc}
G(x,u, \grad u, D^2 u)+\lambda u^{1+\alpha} = f&\ {\rm in}\ \Omega,\\
u=0&{\rm on} \ \partial\Omega.
\end{array}\right.
$$
}

\noindent{\it Furthermore
there exists $\phi>0$ in
$\Omega$ such that $\phi$ is a viscosity solution of 

$$
\left\{\begin{array}{lc}
G(x,\phi, \grad \phi,D^2\phi)+\bar\lambda \phi^{1+\alpha} = 0 & {\rm in}\ 
\Omega\\
\phi=0  & {\rm on}\  \partial \Omega.
\end{array}
\right.
$$
$\phi$ is $\gamma$-H\"older continuous for all $\gamma\in ]0,1[$ 
and locally Lipschitz. }

\bigskip
Let us mention that it is possible to define another "eigenvalue" :
Indeed let 
$$\underline\lambda= \sup\{ \mu, \exists \ \phi<0 \ {\rm in}\ \Omega ,
G(x,\phi,\nabla \phi, D^2\phi)+
\mu |\phi|^\alpha \phi\geq 0\ {\rm in \ the \ viscosity\ sense}\}.$$
If  $\underline G(x,u,p,X) =
-F(x,p, -X)+b(x)\cdot p|p|^{\alpha}+c(x)|u|^\alpha u$ then  $\underline\lambda=\bar  \lambda (\underline G)$. 
 Furthermore if $F$ satisfies (H2) then so does  $\underline F(x,p,X) =
-F(x,p,-X)$.  Hence it is
 possible to prove for $\underline\lambda$ the same results than for
 $\bar\lambda$. It is important to remark that in general 
$\underline G\neq G$ and hence $\underline\lambda$ can be
different from $\bar\lambda$.

\bigskip

While we were completing this paper, we received a very interesting preprint of Ishii and Yoshimura \cite{IY} where similar results are obtained in the case $\alpha=0$. Let us mention that they call the eigenvalue a demi-eigenvalue as in the paper of P.L. Lions \cite{PLL}, and they characterize it  as 
the supremum of those $\lambda\in \R$ for which there is a viscosity supersolution $u\in C(\Omega)$ of $F[u] = \lambda u + 1$ in $\Omega$ which satisfies $u \geq  0$ in $\overline\Omega$. "
(their $F$ is our $-G$). 

\bigskip

In the next section we state precisely the conditions on $G$  and the definition of viscosity solution in this setting. In section 3 we prove a comparison principle and some boundary estimates that allow to prove that for $\lambda<\bar\lambda$ the maximum principle
holds. This will be done in the fourth section, where we also
provide some estimates on $\bar\lambda$ and a further 
comparison principle when
$\lambda<\bar\lambda$. In section 5, using Ishii-Lions technique we prove
regularity results, these in particular give the required relative
compactness for  the sequence of solutions that are used to prove the
main existence's results in the last section.

\section{Main assumptions and definitions. }

In this section, we state the assumptions on the operators
  
$$G(x,u, \nabla u,D^2 u) = F(x,\nabla u, D^2 u)+ b(x).\nabla u |\nabla
u|^\alpha + c(x) |u|^\alpha u $$
treated in this note and the notion of viscosity solution.

The operator $F$ is continuous on $\R^N\times
(\R^N)^\star\times S$, where $S$ denotes the space of symmetric
matrices  on $\R^N$.

The following hypothesis will be considered

\begin{itemize}

\item[(H1)] 
 $F: \Omega\times \R^N\setminus\{0\}\times S\rightarrow\R$, 
and  $\forall t\in \R^\star$,
$\mu\geq 0$,
 $F(x, tp,\mu X)=|t|^{\alpha}\mu F(x, p,X)$.

\item[(H2)]
For $x\in \overline{\Omega}$, $p\in \R^N\backslash \{0\}$, $M\in S$, 
$N\geq 0$
\begin{equation}\label{eqaA}
a|p|^\alpha tr(N)\leq F(x,p,M+N)-F(x,p,M) \leq A
|p|^\alpha tr(N).
\end{equation}

\item[(H3)]
There exists a continuous function $\tilde \omega$, $\tilde \omega(0)=0$
such that for all $x,y,$ $p\neq 0$, $\forall X\in S$
$$|F(x,p,X)-F(y,p,X)|\leq \tilde \omega(|x-y|) |p|^\alpha |X|.$$

\item [(H4)]
 There exists a
continuous function
$  \omega$ with $\omega (0) = 0$, such that if $(X,Y)\in S^2$
and 
$\zeta\in \R$ satisfy
$$-\zeta \left(\begin{array}{cc} I&0\\
0&I
\end{array}
\right)\leq \left(\begin{array}{cc}
X&0\\
0&Y
\end{array}\right)\leq 4\zeta \left( \begin{array}{cc}
I&-I\\
-I&I\end{array}\right)$$
and $I$ is the identity matrix in $\R^N$,
then for all  $(x,y)\in \R^N$, $x\neq y$
$$F(x, \zeta(x-y),p, X)-F(y,  \zeta(x-y),p, -Y)\leq \omega
(\zeta|x-y|^2).$$

\end{itemize}

\bigskip
The condition (H2), usually called uniformly elliptic condition, will be in some cases replaced by the much weaker condition
\begin{itemize}
\item[(H2')]
for all $x\in \Omega$, $p\in \R^N\backslash 0$, $M\in {\mathcal S}$, $N\geq
0$,
$$F(x,p,N+M)\geq F(x,p,M).$$
\end{itemize}

\begin{remark}\label{rem1}

The assumption (H2) and the fact that $F(x,p,0)=0$ implies that 
$$|p|^\alpha (a {\rm tr}(M^+)-A{\rm tr}(M^-))\leq F(x,p,M)\leq |p|^\alpha
(A{\rm tr}(M^+)-a{\rm tr}(M^-))$$
where 
$M=M^+-M^-$ is a minimal decomposition of $M$ into positive and negative
symmetric matrices. 
\end{remark}

We now assume conditions for the lower order terms i.e.
we shall suppose that 
$b:\Omega\mapsto
\R^N$ and $c:\Omega\mapsto\R$ are continuous and bounded.

We shall sometimes require (for example for the comparison and the maximum
principle) that $b$ satisfies:

\begin{itemize}
\item[(H5)]  -Either
$\alpha<0$ and $b$ is H\"olderian of exponent $1+\alpha$,

- or
$\alpha\geq 0$ and, for all $x$ and $y$, 
$$\langle b(x)-b(y), x-y\rangle \leq 0$$
\end{itemize}
\bigskip

Let us recall what we mean by {\it viscosity solutions}, adapted to
our context.

It is well known that in dealing with viscosity respectively  sub and super
solutions one works with  
$$u^\star (x) = \limsup_{y, |y-x|\leq r} u(y)$$

and $$u_\star (x) = \liminf_{y, |y-x|\leq r} u(y).$$
 It is easy to see that 
$u_\star \leq u\leq u^\star$ and $u^\star$ is upper semicontinuous (USC)
$u_\star$ is lower semicontinuous (LSC).  See e.g. \cite{CIL, I}.

\begin{definition}\label{def1}

 Let $\Omega$ be a bounded domain in
$\R^N$, then
$v$,   bounded on $\overline{\Omega}$ is called a viscosity super solution
of
${\mathcal G}(x,\grad u,D^2u)=g(x,u)$ if for all $x_0\in \Omega$, 

-Either there exists an open ball $B(x_0,\delta)$, $\delta>0$  in $\Omega$
on which 
$v= cte= c
$ and 
$0\leq g(x,c)$, for all $x\in B(x_0,\delta)$

-Or
 $\forall \varphi\in {\mathcal C}^2(\Omega)$, such that
$v_\star-\varphi$ has a local minimum on $x_0$ and $\grad\varphi(x_0)\neq
0$, one has
\begin{equation}
{\mathcal G}( x_0,\grad\varphi(x_0),
 D^2\varphi(x_0))\leq g(x_0,v_\star(x_0)).
\label{eq1}\end{equation}

Of course $u$ is a viscosity sub solution if
for all $x_0\in \Omega$,

-Either there exists a ball $B(x_0, \delta)$, $\delta>0$ on which  $u =
cte= c$ and
$0\geq g(x,c)$, for all $x\in B(x_0,\delta)$

-Or  $\forall
\varphi\in {\mathcal C}^2(\Omega)$, such that
$u^\star-\varphi$ has a local maximum on $x_0$ and
$\grad\varphi(x_0)\neq 0$, one has
\begin{equation}
{\mathcal G}( x_0,\nabla \varphi(x_0),
D^2\varphi(x_0))\geq g(x_0,u^\star(x_0)) . \label{eq2}\end{equation}

A	 viscosity solution is a function which is both a super-solution and a sub-solution.
\end{definition}

In particular we shall use this definition with 
${\mathcal G}(x,p,X)=F(x,p,X)+b(x)\cdot p |p|^\alpha $.

See e.g. \cite{CGG} for similar definition of viscosity solution for 
equations with singular operators.

For convenience we recall  the definition of semi-jets given e.g. in 
\cite{CIL}

\begin{eqnarray*}
J^{2,+}u(\bar x) &= &\{ (p, X)\in \R^N\times S, \ u(x)\leq
u(\bar x)+
\langle p, x-\bar x\rangle + \\
& &+{1\over 2} \langle X(x-\bar x), x-\bar
x\rangle
+ o(|x-\bar x|^2) \}
\end{eqnarray*}
and 

\begin{eqnarray*}J^{2,-}u(\bar x)& =& \{ (p, X)\in \R^N\times S, \
u(x)\geq u(\bar x)+
\langle p, x-\bar x\rangle +\\
&&+ {1\over 2} \langle X(x-\bar x), x-\bar
x\rangle
+ o(|x-\bar x|^2\}.
\end{eqnarray*}
In the definition of viscosity solutions the test functions can be
substituted by the elements of the semi-jets in the sense that in the
definition above one can restrict to the functions $\phi$ defined by
$\phi(x)=u(\bar x)+
\langle p, x-\bar x\rangle 
+ {1\over 2} \langle X(x-\bar x), x-\bar
x\rangle$  with $(p,X)\in J^{2,-}u(\bar x)$ when $u$ is a super solution
and $(p,X) \in J^{2,+}u(\bar x)$ when $u$ is a sub solution.  

\bigskip

\begin{remark}\label{dist} In all the paper we shall consider that $\Omega$ is a bounded $\mathcal C^2$ domain. In particular we shall use several times the fact that this implies  that the distance to the boundary:
$$d(x,\partial\Omega):=d(x):=\inf\{|x-y|,\  y\in\partial\Omega\}$$
satisfies the following properties:
\begin{enumerate}
\item $d$ is Lipschitz continuous
\item There exists $\delta>0$ such that in $\Omega_\delta=\{x\in\Omega\ \mbox{such that }\ d(x)\leq \delta\}$, $d$ is $\mathcal C^{1,1}$.
\item $d$ is semi-concave, i.e. there exists $C_1>0$ such that 
$d(x)-C_1|x|^2$ is concave and then $J^{2,+} d(x)\neq \emptyset.$
\item If  $J^{2,-}d(x)\neq \emptyset$, $d$ is  differentiable at $x$
and
$|\nabla d(x)|=1$. 
\end{enumerate}
\end{remark}

\bigskip

\section{A comparison principle   and some boundary estimates.}
 We start by establishing a
comparison result which is a sort of extension of the comparison Theorem
2.1 in
\cite{BD1}.

\begin{theorem}\label{compprnew}
Suppose that $F$ satisfies (H1), (H2'), and (H4),  that $b$ is continuous and
bounded and
$b$ satisfies $(H5)$.  Let $f$ and $g$ be respectively upper and lower semi continuous.

 Suppose that $\beta$ is some continuous function  on
$\R^+$  such that $\beta (0)=0$. Suppose that
$\phi>0$ in $\Omega$  lower semicontinuous and $\sigma$ upper semicontinous, 
 satisfy,  respectively,  in the viscosity sense, 

$$F(x, \grad \phi,D^2\phi)+b(x).\nabla \phi |\nabla
\phi|^{\alpha}-\beta (\phi)\leq f$$
$$F(x,  \grad \sigma,D^2\sigma) +b(x).\nabla
\sigma |\nabla \sigma |^{\alpha} -\beta (\sigma)\geq g.$$ 

Suppose that $\beta$ is increasing on $\R^+$ and $f\leq g$, or $\beta$ is
nondecreasing and $f< g$. 

If
$\sigma\leq\phi$ on
$\partial
\Omega$, then $\sigma \leq \phi$ in $\Omega$. 
\end{theorem}

 Before starting the proof, for convenience of the reader,
let us recall the following lemmata proved in
\cite{BD1}, the first one being  an extension of Ishii's
 acclaimed  result.
\begin{lemma} \label{oldlem1}
Let $\Omega$ be a bounded open set in $\R^N$, which
is piecewise ${\mathcal C}^1$. Let $u$ upper semi-continous in $ \overline{\Omega}$, $v$ lower semicontinous in $ \overline{\Omega}$, $(x_j, y_j)\in \Omega^2$, $x_j\neq y_j$
, and  $q\geq \sup\{2,\frac{\alpha+2}{\alpha+1}\}$.

We assume that the function
$$\psi_j(x,y)=u(x)-v(y)-{j\over q}|x-y|^q$$
 has a local maximum on $(x_j, y_j)$, with $x_j\neq y_j$.
Then, there are
$X_j$,
$Y_j$
$\in {\mathcal S}^N$ such that 
$$ j(|x_j-y_j|^{q-2} (x_j-y_j), X_j)\in
J^{2,+} u(x_j)$$

$$ j(|x_j-y_j|^{q-2} (x_j-y_j), -Y_j)\in J^{2,-}
v(y_j)$$

and $$ -4jk_j\left(\begin{array}{cc} I&0\\
0&I\end{array}\right)\leq \left(\begin{array}{cc} X_j&0\\
0&Y_j\end{array}\right)\leq 3jk_j\left( \begin{array}{cc}
I&-I\\ -I&I
\end{array}
\right)$$ where $$k_j =2^{q-3} q(q-1) |x_j-y_j|^{q-2}
.$$
\end{lemma}

\begin{lemma}\label{dem22}
Under the previous assumptions on $F$, 
let $v$ be a lower semicontinuous, viscosity supersolution of 
$$F(x, \nabla v, D^2
v(x))+b(x).\grad v|\grad v|^{\alpha}-\beta (v(x))\leq f(x)$$ for
some functions
$(f,\beta )$ upper semi continuous in $\Omega$. Suppose that
$\bar x$ is some point in $\Omega$ such that
$$v(x)+C|x-\bar x|^q\geq  v(\bar x),$$
where $\bar x$ is a strict local minimum of the left hand side and $v$ is not locally constant around
$\bar x$. Then, 
$$-\beta( v (\bar x))\leq f(\bar x).$$
\end{lemma}
{\bf Remark:} This Lemma was stated and proved for continuous super solutions in \cite{BD1}. The proof is similar but it is  adjusted  to lower semi continuous super solutions and is given here for the convenience of the reader.

{\bf Proof}
Without loss of generality we can suppose that $\bar x=0$.

Since the infimum is  strict, for $\epsilon>0$, there exists $N$ such
that  for any 
$n>N$ 
 
$$\inf_{{1\over n}\leq |x|\leq R} (v(x)+ C|x|^q)\geq m_n> v(0)+\epsilon$$

We take in the following also $N$ large enough in order that 
$$\left({1\over N}\right)^q C \leq {\epsilon/4}$$
and such that 
$C(diam \Omega+1)^{q-1} {q\over N}< {\epsilon}/4$.

Since $v$ is not locally contant  for any $n$ , there exists
$(t_n, z_n)$ in $B(0,{1\over n})^2$ such that 
$$v(t_n)> v(z_n)+ C|z_n-t_n|^q$$
We consider 
$$\inf_{|x|\leq R} v(x)+ C|x-t_n|^q.$$
We prove in what follows that the infimum is achieved in 
$B(0,{1\over n})$ and that it is  not achieved on $t_n$. 

Let us observe indeed
$$\inf_{|x|\leq {1\over n}} v(x)+ C|x-t_n|^q\leq v(0)+ \epsilon/4$$

and since 
$$\inf_{|x|\leq {1\over n}} v(x)+ C|x-t_n|^q\leq v(z_n)+ C|z_n-t_n|^q<
v(t_n)$$
the infimum on $|x|\leq 1/n$ cannot be achieved on $t_n$. 

Moreover 
\begin{eqnarray*}
\inf_{{1\over n}\leq |x|\leq R} (v(x)+ C|x-t_n|^q)&\geq&
\inf_{{1\over n}\leq |x|\leq R} (v(x)+ C|x|^q+ C|x-t_n|^q-C|x|^q)\\
&\geq&
m_n-Cq|t_n| |x-\theta t_n|^{q-1} \geq m_n-C(diam \Omega+1)^{q-1}|t_n|\\
&\geq&
m_n-\epsilon/4\geq v(0)+{3\epsilon\over 4}
\end{eqnarray*}

Since the infimum cannot be achieved on $t_n$,  let $y_n$, $|y_n|\leq
{1\over n}$  be a point such that the infimum is achieved on  $y_n$, then

$$v(x)+ C|x-t_n|^q\geq v(y_n)+ C|y_n-t_n|^q$$
and then 
$$\varphi(x) = v(y_n)+ C|y_n-t_n|^q-C |x-t_n|^q$$
is a test function for $v$ on $y_n$ with a gradient  $\neq 0$ on that
point. 
Since $v$ is a supersolution one gets 
$$-A C^{\alpha+1} |y_n-t_n|^{q(\alpha+1)-\alpha-2}-\beta(v(y_n))\leq
f(y_n)$$

Let us observe that $v(y_n)\rightarrow v(0)$. Indeed one has by the lower
semicontinuity of $v$ 
$$v(0)\leq \liminf v(y_n)$$
and using 
$$v(0)+ C|t_n|^q\geq v(y_n)+ C|y_n-t_n|^q$$
one has the reverse inequality. 

Then by  the uppersemicontinuity of $f$ and $\beta$ one gets
that 
$$-\beta(v(0))\leq f(0)$$
which is the desired conclusion.

\bigskip
\noindent {\bf Proof of Theorem \ref{compprnew}.}

Suppose by contradiction that max $(\sigma-\phi)>0$ in $\Omega $. 
Since $\sigma\leq \phi$ on the boundary, the supremum can only be
achieved inside $\Omega$. 

Let us consider for
$j\in
\N$ and for some
$q>
\max ( 2, {\alpha+2\over
\alpha+1})$  $$\psi_j (x,y)
= \sigma(x)-\phi(y)-{j\over q} |x-y|^q.$$

\bigskip
 Suppose that $(x_j, y_j)$
is a maximum for $\psi_j$.  Then

(i) from the boundedness of $\sigma$ and $\phi$ one deduces that $|x_j-y_j |\rightarrow 0$ as $j\rightarrow \infty$. Thus up to subsequence $(x_j, y_j)\rightarrow (\bar x,\bar x)$

 (ii) One has $\lim\inf \psi_j(x_j,y_j)\geq \sup  (\sigma-\phi)$;
 
  (iii) $\lim\sup \psi_j(x_j,y_j)\leq\lim\sup \sigma(x_j)-\phi(y_j)=  \sigma(\bar x)-\phi(\bar x) $
  
 (iv) Thus $j|x_j-y_j|^q\rightarrow 0$ as $j\rightarrow +\infty$ and $\bar x$ is a maximum point for $\sigma-\phi$.

\bigskip
{\bf Claim: } {\it For $j$ large enough,
 there exist $x_j$ and $y_j$ such that $(x_j, y_j)$ is a maximum pair
for $\psi_j$ and 
$x_j\neq y_j$.}

 \noindent Indeed suppose that $x_j = y_j$. Then one would have
\begin{eqnarray*}
\psi_j (x_j, x_j)&=& \sigma(x_j)-\phi(x_j)\\
&\geq& \sigma(x_j)-\phi(y)-{j\over q} |x_j-y|^q;\\
\psi_j (x_j, x_j)&=& \sigma(x_j)-\phi(x_j)\\
&\geq& \sigma(x)-\phi(x_j)-{j\over q} |x-x_j|^q;\\
\end{eqnarray*}
 and then $x_j$ would be a local maximum for
$$\Phi:=\phi(y)+{j\over q} |x_j-y|^q.$$ 
and similarly  a local minimum  for
$$\Sigma:=\sigma(x)-{j\over q} |x_j-x|^q.$$

We first exclude that $x_j$ is  both a strict local maximum and a strict local minimum.
Indeed  in that case, by Lemma \ref{dem22}
$$-\beta(\phi(x_j))\leq f(x_j)$$

$$-\beta(\sigma (x_j))\geq g( x_j)$$

This is a contradiction because either $\beta$ is increasing
$$-g(x_j)\geq \beta(\sigma (x_j))> \beta(\phi(x_j)\geq -f( x_j)\geq
-g(x_j)$$ or
$$-g(x_j)\geq \beta(\sigma (x_j))\geq  \beta(\phi(x_j))\geq -f( x_j)> -g(x_j).$$
Hence  $x_j$ cannot be   both a strict minimum for $\Phi$  and a 
strict maximum for $\Sigma$.
 In the first case  there exist
$\delta>0$ and $R> \delta$ such that $B(x_j, R)\subset \Omega$ and 

$$\phi(x_j) = \inf_{\delta\leq |x-x_j|\leq R} \{ \phi(x)+ {j\over q}
|x-x_j|^q\}.$$

Then if $y_j$ is a point on which the minimum above is achieved, one has 
$$ \phi(x_j)=\phi(y_j)+{j\over q} |x_j-y_j|^q,$$ 
and  $(x_j, y_j)$ is still a maximum point for $\psi_j$ since for all
$(x,y)\in \Omega^2$

$$\sigma(x_j)-\phi(y_j)-{j\over q} |x_j-y_j|^q = \sigma(x_j)-\phi(x
_j)\geq
\sigma (x)-\phi(y)-{j\over q} |x-y|^q.$$
This concludes the Claim.
In the other case, similarly, one can replace $x_j$ by a point $y_j$ near $x_j$ with  
$$ \sigma (x_j)=\sigma (y_j)+{j\over q} |x_j-y_j|^q,$$ 
and  $(y_j, x_j)$ is still a maximum point for $\psi_j$.

\bigskip We can now conclude. By Lemma \ref{oldlem1}   there exist
$X_j$ and $Y_j$ such that 

$$\left(j|x_j-y_j|^{q-2} (x_j-y_j), {X_j}\right)\in J^{2,+} \sigma(x_j)$$
and 

$$\left(j|x_j-y_j|^{q-2} (x_j-y_j), {-Y_j}\right)\in J^{2,-} \phi(y_j).$$
We can use the fact that $\sigma$ and $\phi$ are respectively sub and super solution to obtain:

\begin{eqnarray*}
g(y_j)&\leq & F(y_j, j|x_j-y_j|^{q-2} (y_j-x_j), -Y_j)+\\
&+ & b(y_j). j^{1+\alpha}|x_j-y_j|^{(q-1)(1+\alpha)-1}
(y_j-x_j)-\beta(\phi(y_j))\\
 &\leq &F(x_j, j|x_j-y_j|^{q-2} (x_j-y_j),X_j)+ \\
&+& b(x_j).j^{1+\alpha}|x_j-y_j|^{(q-1)(1+\alpha)-1}
(x_j-y_j)+\\
&+&  C(j|x_j-y_j|^q)^{1+\alpha}+\omega (j|x_j-y_j|^{q}+
{1\over j})-\beta (\phi(y_j))\\ 
&\leq & f(x_j) + C(j|x_j-y_j|^q)^{1+\alpha}+\omega (j|x_j-y_j|^{q}+
{1\over j}) +\beta(\sigma(x_j))-\beta(\phi(y_j)).
\end{eqnarray*}

Passing to the limit and using the fact that $g$ and $f$ 
are respectively lower and upper semi continuous and $\beta$ is
continuous, we obtain 

$$g(\bar x)\leq f(\bar x)+\beta(\sigma(\bar x))-\beta(\phi(\bar x))$$
which  contradicts our hypotheses in all cases and $\sigma\leq 0$ in $\Omega$.
This ends the proof.

\bigskip

As an application of  the comparison theorem
 we will state bounds for sub and super solutions near the boundary.
The conclusions given in  Proposition \ref{ah} and Corollary \ref{minsup} will
be used in the proof of the maximum principle Theorem \ref{maxp}.

\begin{proposition} \label{ah} Suppose that $F$ satisfies (H1) and (H2), and
that 
$b$ is bounded.

Let $u$ be  uppersemicontinuous  subsolution of 

$$\left\{
\begin{array}{lc}
F(x, \grad u,D^2u)+b(x).\grad u |\grad u|^{\alpha}\geq
 -m &
{\rm in}\ \Omega\\
u =0 & {\rm on}\ \partial\Omega
\end{array}
\right.
$$ 
for some constant $m\geq 0$. 
Then there exists $\delta>0$ and some constant $C_3$  that depends only on the structural
data such that
 $$u(x)\leq C_3d(x,\partial \Omega)$$
if  the distance to the boundary $d(x,\partial\Omega)< \delta$.

\end{proposition}
{\bf Proof}  Let $d(x)=d(x,\partial\Omega)$.
First let us observe that one can assume that there exists $d_0$ such
that  $\Omega_{d_0}=\{x\in\Omega\quad \mbox{such that}\quad d(x)< d_0\}$ the supremum of 
$u$ is positive  because otherwise there is nothing to prove.

We recall that  $\Omega$ is a bounded $C^2$ domain  and using the properties of the distance function
stated in Remark \ref{dist} we know that 
 $D^2d\leq C_1 Id$ for some constant $C_1$ that depends on $\Omega$.
Let $\Omega_\delta=\{x\in\Omega; \ d(x)<\delta\}$. Suppose that $\delta<
{a\over 4(C_1(A+a)N+|b|_\infty)}$.  For some constants
$\gamma$ and $C_2$ that will be chosen later we introduce
$$\psi(x)=\gamma \log (1+C_2d(x)).$$

We use the inequalities 
$$tr(D^2 \psi)^+\leq {\gamma C_2 N C_1\over 1+C_2 d}$$
and 
$$tr(D^2 \psi)^-\geq {\gamma C_2^2\over (1+C_2d)^2}-{\gamma C_2
C_1 N\over 1+C_2 d}$$

 Then choosing $C_2>4\frac{C_1(A+a)N+|b|_\infty}{a}$, one has

\begin{eqnarray*}
& & F(x,\grad\psi,D^2\psi)  +  b(x).\grad\psi|\grad \psi|^\alpha \\
 &\leq&
\gamma^{\alpha+1}\left(\frac{C_2}{1+C_2d}\right)^{\alpha+1}
\left[\frac{-aC_2}{1+C_2d}+C_1(A+a)N+|b|_\infty\right]\\
 & \leq&   - \frac{a}{2} \gamma^{\alpha+1}\left(\frac{C_2}{1+C_2\delta}\right)^{\alpha+2}.
\end{eqnarray*}

Now we choose $\gamma$ sufficiently large that $\gamma> \sup_{\{x, \ d(x)\leq
\delta\}} u(x)$ and 

$$ \frac{a}{4}\gamma^{\alpha+1}\left(\frac{C_2}{1+C_2\delta}\right)
^{\alpha+2}\geq m$$ 
i.e.

$$\gamma=\max\left\{(\frac {am}{ 4})^{\frac{1}{\alpha+1}}\left(\frac{1+C_2\delta}{C_2}\right)
^{\frac{\alpha+2}{\alpha+1}}
,\ \ \sup u_{\{x, \ d(x)\leq\delta\}} \ \right\}.$$

With this choice of constants we have obtained that 
\begin{eqnarray*}
& & F(x,\grad\psi,D^2\psi)  +  b(x).\grad\psi|\grad \psi|^\alpha \\
&\leq &- \frac{a}{2}\gamma^{\alpha+1}\left(\frac{C_2}{1+C_2\delta}\right)^{\alpha+2}\\
&< & - \frac{a}{4}\gamma^{\alpha+1}\left(\frac{C_2}{1+C_2\delta}\right)^{\alpha+2}\\
&\leq& -2m\\
&\leq&  F(x,\grad u,D^2 u) +b(x)\cdot\grad u|\grad u|^{\alpha+1}
\end{eqnarray*}
and furthermore $u\leq \psi$ on $\partial\Omega_\delta$. 

Hence by Theorem \ref{compprnew} with 
$f=- \frac{a}{2}\gamma^{\alpha+1}\left(\frac{C_2}
{1+C_2\delta}\right)^{\alpha+2}$, $g(x)=
-m$ we obtain 
$$u(x)\leq \gamma  \log (1+C_2 d(x))\leq \gamma C_2 d(x)$$
since  $u\leq\psi$ in $\Omega_\delta$.
This ends the proof.

\bigskip

The comparison principle in \cite{BD1} allows also to establish a strict maximum 
principle:

\begin{theorem}\label{smax}

Suppose that $F$ satisfies (H2), $b$ and $c$  are continuous and 
bounded and $b$ satisfies $(H5)$. Let
$u$ be a viscosity non-negative lowersemicontinuous  super
solution of 
$$F(x,\nabla u, D^2 u)+ b(x).\nabla u |\nabla u|^\alpha +c(x) u^{\alpha+1}
\leq 0.$$

 Then either $u\equiv 0$ or $u>0$ in $\Omega$. 
\end{theorem}
{\bf Remark:} Other strong maximum principles  and strong minimal principles have been established in \cite{BDL} for a more general class of fully nonlinear operators that are "proper" .

{\bf Proof.}
Using the inequality in (H2), let us
recall,  using Remark \ref{rem1}, that
\begin{eqnarray*}
F(x,  p, M) &\geq& |p|^\alpha (a tr (M)^+-A
 tr (M)^-)\\
&:=& H(p, M).
\end{eqnarray*}
Hence it is sufficient to prove the proposition when $u$ is
 a super solution of
$$H(\nabla u, D^2 u)+b(x).\nabla u|\nabla u|^\alpha + c(x)u^{1+\alpha}
=0.$$

$H$ does not depend on
$x$ and it satisfies the hypothesis of Theorem \ref{compprnew}.

Moreover one can assume that $c$ is some negative constant.
Indeed, suppose that we have proved that for any $u\geq 0$ super solution of
\begin{equation}\label{306}
H(\nabla u, D^2 u)+b(x).\nabla u |\nabla u|^\alpha-|c|_\infty  u^{\alpha+1} \leq 0,
\end{equation}
we have that $u>0$ in $\Omega$.  Then  if 
$$F(x, \nabla u, D^2 u)+ b(x).\nabla u|\nabla u|^\alpha +
c(x)u^{\alpha+1}
\leq 0$$
for some $u\geq 0$ we have that $u$ is a non negative super solution of (\ref{306}) and then 
$u>0$ and that would conclude the proof.

Hence we suppose by contradiction that  $x_0$ is some point inside
$\Omega$ on which $u(x_0)=0$. Following e.g. Vazquez \cite{V}, one can
assume that on the ball $|x-x_1| = |x-x_0| = R$, $x_0$ is the only point
on which $u$ is zero and that $B(x_1, {3R\over 2})\subset \Omega$.
Let $u_1 = \displaystyle\inf_{|x-x_1| = {R\over 2}} u>0$, by the
lower semicontinuity of
$u$. Let us construct a sub solution on the annulus ${R\over
2}\leq |x-x_1|=\rho<{3R\over 2}$.

Let us recall that if $\phi (\rho) =
e^{-k\rho}$, the eigenvalues of $D^2\phi$ are
$\phi^{\prime\prime}(\rho)$ of multiplicity 1 and $ \displaystyle{\phi^\prime\over
\rho}$ of multiplicity N-1.

 Then take  $k$ such that
$$k^{\alpha+2}>\left({2(N-1)A\over R a}+ |b|_\infty\right)k^{\alpha+1}
+ |c|_\infty.$$

 If $k$ is as above, let  $m$ be  chosen
such that $$m(e^{-kR/2}-e^{-kR})= u_1$$  and define  $v(x) =m
(e^{-k\rho}-e^{-kR})$ with $\rho=|x|$. The function  $v$ is a  strict 
subsolution in the annulus, in the sense that it satisfies  
$H(\nabla v, D^2v) +b(x).\nabla v |\nabla v|^{\alpha}-|c|_\infty
v^{\alpha+1}>0$ in the annulus. Furthermore

$$\left\{\begin{array}{cc}
v\leq u &{\rm on} \ |x-x_1|=\displaystyle{R\over 2}\\
v\leq 0\leq u&{\rm on}\  |x-x_1|=\displaystyle{3R\over 2}.
\end{array}\right.$$

 Hence $u\geq v$ everywhere on the boundary of the annulus. In fact $u\geq v$ everywhere in the annulus, since we can use  the comparison principle Theorem \ref{compprnew} for the operator
$H+b(x).\nabla . |\nabla .|^{\alpha}-|c|_\infty |.|^{\alpha}$. 
Then
$v$ is a test function for $u$ at $x_0$.  Then, since $u$ is a super solution and
$\nabla v(x_0)\neq 0$:
$$H( \nabla v(x_0), D^2v(x_0))+b(x_0).\nabla v(x_o) |\nabla v(x_o)|^{\alpha}-|c|_\infty
v^{\alpha+1}(x_0)\leq  0$$
which clearly contradicts the definition of $v$. Finally $u$ cannot be zero
inside $\Omega$. This ends the proof.

\begin{corollary}[Hopf] \label{Hopf}Let
$v$ be a viscosity continuous super solution of 
$$F(x,\nabla v, D^2 v)+ b(x).\nabla v |\nabla v|^\alpha +c(x) |v|^{\alpha}v
\leq 0.$$
Suppose that $v$ is positive in a neighborhood of $x_o\in\partial\Omega$ and
$v(x_o)=0$ then there exist $C>0$  and $\delta>0$ such that
$$v(x)\geq C|x-x_o|$$

for $|x-x_o|\leq \delta$.
\end{corollary}

To prove this corollary just proceed as in the proof of Theorem 
\ref{smax} and remark that $e^{-k\rho}-e^{-kR}\geq C(R-\rho)$.

In fact, one can get a better estimate about supersolutions near the
boundary i.e. some sort of limited expansion at the order two. We still denote by  $d(x)$ the distance to the boundary of $\Omega$ and, for $d>0$,  $\Omega_d=\{x\in\Omega,\quad d(x)\leq d\}$.

\begin{proposition}\label{minsup}
Suppose that $v$ is a lowersemicontinuous supersolution  
of 
$$F(x,\nabla v,D^2 v ) +b.\nabla v|\nabla
v|^\alpha +c|v|^{\alpha}v\leq 0$$
in $\Omega_{d_1}$, for some $d_1>0$,   
and $v$ is $\geq 0$ on the boundary, $v>0$. Then there exists $d_2\leq d_1$, $d_2>0$   and
some constants
$\gamma $, $C>0$ such that on  $\Omega_{d_2}$
$$v(x)\geq \gamma (d(x)+ \log (1+Cd(x)^2))\geq \gamma (d(x)+{Cd(x)^2\over 2}).$$
\end{proposition}
{\bf Proof.} We start by proving the following

{\bf Claim}: { For some constant } $C>0$  large enough, there exists a neighborhood of
$\partial \Omega$ such that 
$$\varphi(x) = d(x)+\log (1+Cd^2(x))$$
satisfies 
$$F(x,\nabla\varphi,D^2 \sigma ) +b.\nabla \varphi |\nabla
\varphi|^\alpha +c|\varphi|^{\alpha}\varphi>m>0$$
for some constant $m>0$. 

\bigskip
Let $d_0$ be such that in $\Omega_{d_0}:=\{x\in\Omega :\quad d(x)< d_0\}$ 
the distance is smooth and there exists $C_1$ such that $|D^2d|_\infty\leq C_1$ as seen in Remark \ref{dist}. Note that this implies that
$tr(D^2d)^++ tr (D^2d)^-\leq C_1N$.  

Let $C> ({25\over 6} +{5A\over a} ) NC_1+ 25 {|b|_\infty 2^{1+\alpha} +
|c|_\infty\over
\inf (1,2^\alpha)}$
 and $d< \inf ({1\over 2C}, {1\over 2}, d_0)$.

In $\Omega_{d_0}$
$$|\varphi|\leq d+ Cd^2 \leq {1\over 2}+{1\over 4}\leq 1$$
 We
compute the two first derivatives of $\varphi$  :
$$\nabla \varphi = \nabla d (1+{2Cd\over 1+Cd^2})$$
and then 
$$1\leq |\nabla \varphi|\leq 2$$
$$D^2 \varphi = D^2 d (1+{2Cd\over 1+Cd^2}) +
{2C(\nabla d\otimes \nabla d)(1-Cd^2)\over(1+Cd^2)^2.}$$
In particular
$$(D^2 \varphi)^- \leq  (D^2d)^- (1+{2Cd\over
1+Cd^2}),$$
and
$$(D^2\varphi)^+ \geq  2C(\nabla d\otimes\nabla d){(1-Cd^2)\over
(1+Cd^2)^2} -(D^2 d)^-(1+{2Cd\over
1+Cd^2}),.$$
These imply that 
$$tr(D^2\varphi)^-\leq 2C_1N$$
and 
$$tr (D^2\varphi)^+\geq C{24\over 25}-2NC_1\geq {12C\over 25}.$$

Hence we obtain
\begin{eqnarray*}
F(x, \nabla \varphi, D^2 \varphi) &+ & b(x).\nabla
\varphi|\nabla \varphi|^\alpha + c(x)\varphi^{1+\alpha}\\
&\geq&
\inf (1 ,2^\alpha)({12 aC\over 25}-2A NC_1)-|b|_\infty |\nabla
\varphi|^{1+\alpha}-|c|_\infty |\varphi|^{1+\alpha} \\
&  \geq& \inf
(1,2^\alpha) {2aC\over 25} -|b|_\infty 2^{1+\alpha} -|c|_\infty \\
&\geq &
\inf (1, 2^\alpha) ({aC\over 25})
>0
\end{eqnarray*}
This ends the proof of the Claim.

To conclude the proof of the proposition
 we choose $C$ and $d_0$ as in the claim, $d_2\leq
(d_1, d_0)$ . Since 
$v>0$ inside
$\{x,\ d(x)< d_1\}$  let $\gamma $ be such that $\gamma (d_2+\log
(1+Cd_2^2))\leq min_{d(x) = d_2} v$. Then 
$v\geq \gamma(d+\log (1+Cd^2))$ on the boundary of the "crown " $\{x, 0<
d(x)<d_2\}$ in $\Omega$. Since in addition $v$ satisfies 
$$F(x,\nabla v,D^2 v ) +b.\nabla v|\nabla
v|^\alpha +c|v|^{\alpha}v\leq 0$$
and $\varphi= \gamma (d+\log (1+Cd^2))$ satisfies 
$$F(x,\nabla \varphi,D^2 \varphi ) +b.\nabla \varphi|\nabla
\varphi|^\alpha +c|\varphi|^{\alpha}\varphi> 0,$$
the comparison principle implies that $v\geq \gamma (d+\log (1+Cd^2))$ in
$\{x, d(x)\leq d_2\}$. 
This ends the proof.

\section { Maximum principle  for $\lambda<\bar\lambda$;  bounds for
$\bar\lambda$.}

\subsection{Maximum principle.}
We can now state and prove the following Maximum principle:

\begin{theorem}\label{maxp}

Let $\Omega$ be a bounded domain of $\R^n$.
Suppose that $F$ satisfies (H1), (H2'), (H4),  that $b$ and $c$ are continuous and
$b$ satisfies (H5).  Suppose that
$\tau<\bar\lambda$ and that
$u$ is a  viscosity sub solution of 

$$G(x, u, \nabla u, D^2 u)+ \tau |u|^{\alpha}u \geq 0\  \ {\rm in}\ \ \Omega$$
with $u\leq 0$ on the boundary of $\Omega$, then 
$u\leq 0$ in $\Omega$. 

\end{theorem}
\noindent {\bf Remark:} Similarly it is possible to prove that if $\tau<
\underline\lambda$ and  $v$ is a super solution of 
$$G(x, v, \nabla v, D^2 v)+ \tau |v|^{\alpha}v \leq 0\  \ {\rm in}\ \ \Omega$$
with $v\geq 0$ on the boundary of $\Omega$ then $v\geq 0$ in $\Omega$.
\bigskip

{\bf Proof.}
Let $\lambda\in ]\tau, \bar\lambda[ $, and let 
$v$ be a   super solution of 
$$G(x,v,\nabla v, D^2v) +\lambda v^{\alpha+1}\leq 0,$$
satisfying $v>0$ in $\Omega$, which exists by definition of $\bar\lambda$.

We assume by contradiction that $\sup u(x)>0$ in $\Omega$. 

{\bf Claim} :  $\sup\frac{u}{v}<+\infty$.

Near the boundary this holds true since from Proposition \ref{ah} and Corollary \ref{Hopf} 
there exists $\delta>0$ such that $u(x)\leq Cd(x)$ and $v(x)\geq
C^\prime d(x)$ for $d(x)\leq\delta$, for some constants $C$ and
$C^\prime$.  

In the interior we just use the fact that 
$v\geq C^\prime\delta>0$ in  $\Omega_\delta=\{x:\ d(x)\geq \delta\}$ and we can conclude that 
 ${u\over v}$  is
bounded in $\Omega$.


\bigskip

We now define $\gamma^\prime =\sup_{x\in {\Omega}}{u\over v}$
and
$w=\gamma v$, where $0<\gamma<\gamma^\prime$ and $\gamma$ is 
sufficiently close to $\gamma^\prime$ 
in order that 
$\displaystyle{\lambda-\tau \left({\gamma^\prime\over \gamma}\right)^{1+\alpha}\over
\left({\gamma^\prime\over \gamma}\right)^{1+\alpha}-1}\geq 2|c|_\infty$.
Furthermore by definition of the supremum there exists $\bar
y\in\bar \Omega$ such that
$\sup \frac{u(x)}{v(x)}= \frac{u(\bar y)}{v(\bar y)}= \gamma^\prime $.

Clearly, by homogeneity, $G(x,w,\grad w,D^2 w)+\lambda w^{1+\alpha}\leq 0$. 

The supremum of $u-w$ is strictly positive, and it is necessarily 
achieved on $\bar x \in \Omega$ since on the boundary $u-w\leq 0$.
One has 

$$(u-w)(\bar x)\geq (u-w)(\bar y)$$
and then 
$$w(\bar x)\leq u(\bar x)+ (w-u)(\bar y) < u(\bar x).$$
 On the other hand 

$$ (u-w)(\bar x)\leq (\gamma^\prime-\gamma) v(\bar x)$$
and then 
$$w(\bar x)\geq {\gamma \over \gamma^\prime} u(\bar x).$$

As in the comparison principle, we consider, for $j\in \N$ and for some $q> \max ( 2, {\alpha+2\over
\alpha+1})$: 
 $$\psi_j (x,y)
= u(x)-w(y)-{j\over q} |x-y|^q.$$  
Since $\sup (u-w)>0$,  the supremum of $\psi_j$ is achieved in $(x_j,y_j)\in\Omega^2$.

For $j$ large enough,
$\psi_j$ achieves its positive maximum on some couple
$(x_j, y_j)\in
\Omega^2$ such that 

1) $x_j\neq y_j$ for $j$ large enough, (this uses lemma \ref{dem22} and the definition of $\gamma$).

2) $(x_j, y_j)\rightarrow (\bar x,\bar x)$ which is a maximum point for ${u - w}$ and it is an interior point

3) $j|x_j-y_j|^q\rightarrow 0$,

4) there exist $X_j$ and $Y_j$ in $S$ such that 
$$\left(j|x_j-y_j|^{q-2} (x_j-y_j), X_j\right)\in J^{2,+} u(x_j)$$
and 
$$\left(j|x_j-y_j|^{q-2} (x_j-y_j), -Y_j\right)\in J^{2,-} w(y_j)$$
Furthermore
$$\left(\begin{array}{cc}
X_j&0\\
0&Y_j
\end{array}\right)\leq j\left(\begin{array}{cc}
D_j&-D_j\\
-D_j&D_j
\end{array}\right)\leq 2^{q-2} j q(q-1) |x_j-y_j|^{q-2}
\left(\begin{array}{cc}
I&-I\\
-I&I\end{array}\right)
$$ with 
$$D_j = 2^{q-3} q|x_j-y_j|^{q-2} (I+{(q-2)\over |x_j-y_j|^2}
(x_j-y_j)\otimes (x_j-y_j)).$$ 
The proof of these facts proceeds similarly to the one given in  Theorem \ref{compprnew}.

\bigskip

Condition (H4) implies that

$$F(x_j, j(x_j-y_j)|x_j-y_j|^{q-2} , {X_j})-F(y_j, 
j(x_j-y_j)|x_j-y_j|^{q-2}
, {-Y_j})\leq
\omega (j|x_j-y_j|^q).$$

 Then, using the above inequality, the properties of the 
sequence $(x_j,y_j)$, the condition on $b$ -with $C_b$ below being either
its  H\"older constant or $0$-, and the homogeneity condition (H1), one
obtains 

\begin{eqnarray*} 
-(\tau +c(x_j))u (x_j)^{1+\alpha} &\leq& F(x_j,j(x_j-y_j)|x_j-y_j|^{q-2} , X_j)\\
& &+b(x_j).j^{(1+\alpha)}|x_j-y_j|^{(q-1)(1+\alpha)-1} (x_j-y_j)\\
 &\leq &F(y_j, 
j(x_j-y_j)|x_j-y_j|^{q-2}, -Y_j)\\
&& +{\omega (j|x_j-y_j|^q)} + C_bj^{1+\alpha}|x_j-y_j|^{q(1+\alpha)}\\
&&+b(y_j).j^{1+\alpha}|x_j-y_j|^{(q-1)(1+\alpha)-1} (x_j-y_j)+o(1)\\
&\leq &(-\lambda-c(y_j)) w(y_j)^{1+\alpha}+o(1).
\end{eqnarray*}

By passing to the limit when $j$ goes to infinity, since $c$ is continuous one gets

$$-(\tau+c(\bar x)) u(\bar x)^{1+\alpha}\leq-(\lambda 
+c(\bar x))w(\bar x)^{1+\alpha}$$

If $c(\bar x)+\lambda>0$ one obtains that 
$$-(\tau+c(\bar x) ) u(\bar x)^{1+\alpha} \leq -(\lambda + c(\bar x))
\left({\gamma \over \gamma^\prime}\right)^{1+\alpha} u(\bar x)^{1+\alpha}$$
 This
contradicts the hypothesis that
${\lambda-\tau \left({\gamma^\prime\over \gamma}\right)^{1+\alpha}\over \left({\gamma^\prime\over
\gamma}\right)^{1+\alpha}-1}\geq 2|c|_\infty$. 

If $c(\bar x)+\lambda = 0$ then $\tau< \lambda$ implies that 
$$0<-(\tau+c(\bar x)) u(\bar x)^{1+\alpha} \leq 0$$
a contradiction. Finally if  $c(\bar x)+\lambda <0$ we obtain 
$$-(\tau+c(\bar x)) u(\bar x)^{1+\alpha}\leq -(\lambda+ c(\bar x)) 
w(\bar x)^{1+\alpha}\leq-(\lambda+ c(\bar x)) u(\bar x)^{1+\alpha}$$
once more a contradiction since $\tau< \lambda$. 
This ends the proof.

\subsection{Bounds on $\bar\lambda$}

\begin{proposition}\label{propmin} Let $c(x)\equiv 0$.
Let $F$ satisfying (H1), (H2). Suppose that  $\Omega$ is bounded in at
least one  direction, say $e_1$,   i.e. there exists $R$  such that $\Omega
\subset [-R,R]\times \R^{N-1}$, and that $b_1(x)=<b(x),e_1>$  is
bounded.  Then there exist some constants $C_1>0$, $C_2>0$ which depend on
$a$ and $N$ such that 
$$\bar\lambda> {C_1e^{-C_2|b_1|_\infty}\over R^{2+\alpha}}.$$
\end{proposition}

We deal with the particular case of the dimension 1. In that case we
shall use variational techniques and weak solutions to estimate the first
eigenvalue, this  being justified by the following lemma 

\begin{lemma}

Suppose that $\Omega= ]-R,R[,$ $R>0$, that $a$ is some continuous
function such that 
$$0<a\leq a(x)\leq A$$
in $[-R,R]$, and that $b$ is continuous and bounded. 
Suppose that $g$ is continuous, then the weak solutions (in
$ W^{1,2+\alpha} (]-R,R[)$) and the viscosity solutions of 
$$a(x) |u^\prime |^\alpha u^" + b(x) |u^\prime|^\alpha u^\prime = g(x)$$
for $x\in ]-R,R[$
coincide.
\end{lemma}
{\bf Proof.} 
Suppose first that $u\in W^{1,2+\alpha}$ is a weak solution.

  The previous
equation can also be written as 
$${d\over dx} (|u^\prime |^\alpha u^\prime e^{\int_0^x{(\alpha+1)b(t)\over a(t)}
dt} ) = g(x) e^{\int_0^x{(\alpha+1)b(t)\over a(t)}
dt} $$
Then since $|u^\prime|^\alpha u^\prime= h\in L^{\alpha+2\over
\alpha+1}$ and $e^{\int_0^x{b(t)\over a(t)}
dt}$ is continuous, the product is a distribution $T$  which satisfies 
"$T^\prime $ is continuous". Then $T$ is ${\mathcal C}^1$, hence  $h(x) = Te^{-\int_0^x{b(t)\over a(t)} dt}$ is ${\mathcal C}^1$. Finally $u^\prime (x) =
h(x)^{1\over 1+\alpha}$ is  ${\mathcal C}^1$ on every point where $h(x)\neq
0$, i.e. on each point where $u^\prime(x)\neq 0$.  Finally $u$ is
${\mathcal C}^2$ on such point, and then on those points it satisfies the equation in the classical sense. 

We now prove that $u$ is a viscosity solution. 

For that aim let $\varphi$ be such that $(u-\varphi)(x)\geq 0 =
(u-\varphi)(\bar x)$ for all $x$ in a neighborhood of $\bar x$. Since $u\in {\mathcal C^1}$, $\varphi'(\bar x)=u'(\bar x)$. If $\varphi^\prime (\bar x)= 0$ then there is nothing to test, if 
$\varphi^\prime (\bar x)\neq 0$ then $u^{\prime\prime}(\bar x)$
exists. Moreover $\varphi^{\prime\prime}(\bar x)\leq u^{\prime\prime}(\bar x)$, and then 

$$a(\bar x)|\varphi^\prime |^\alpha \varphi^{\prime\prime}(\bar x) + b(x)
|\varphi^\prime (\bar x)|^\alpha \varphi^\prime (\bar x)\leq a(\bar
x)|u^\prime  |^\alpha u^{\prime\prime}(\bar x) + b(x) |u^\prime (\bar
x)|^\alpha u^\prime (\bar x)\leq g(\bar x)$$
one sees that $u$ is a super-solution. 

Suppose that $\varphi$ is some test function  by above  for $u$ on $\bar
x$, again we are only intrested in the case $\varphi^\prime (\bar x)\neq 0$ which  implies that $u^\prime $ cannot be
zero, and $\varphi^{\prime\prime}(\bar x)\geq u^{\prime\prime}(x)$
Then since  on those points $u$  is a solution in the classical
sense 

$$a(\bar x)|\varphi^\prime|^\alpha \varphi^{\prime\prime}(\bar x) + b(x)|\varphi^\prime
|^\alpha \varphi^\prime (\bar x)\geq g(\bar x).$$
This implies 
 which implies that $u$ is a sub solution .

We prove that the viscosity solutions are weak solutions, in the one
dimensional case. 

Let $v$ be a weak solution of 
$$a(x)|v^\prime|^\alpha  v^" + b(x) |v^\prime|^\alpha v^\prime = g(x)$$
$v= 0$ on the boundary. 

Let now $u$ be a viscosity solution of the same equation. We want to
prove that $u=v$. 
For that aim let $\epsilon$ and let $v_\epsilon$ be the weak solution of 
$$a(x)|v_\epsilon^\prime|^\alpha  v_\epsilon^" + b(x)
|v_\epsilon^\prime|^\alpha v_\epsilon^\prime = -\epsilon + g(x),$$
$v_\epsilon= 0$ on the boundary,
and $v^\epsilon$ be the weak solution of 
$$a(x)|(v^\epsilon)^\prime|^\alpha  (v^\epsilon)^" + b(x)
|(v^\epsilon)^\prime|^\alpha (v^\epsilon)^\prime =
\epsilon+ g(x),$$
$v^\epsilon= 0$ on the boundary. 
By the previous part $v_\epsilon$ and $v^\epsilon$ are viscosity
solutions and by the comparison theorem \ref{compprnew} gets 
$$v^\epsilon\leq u\leq v_\epsilon.$$
Moroever by passing to the limit  for weak solutions (for example using
variational technics) it is easy to prove that $v_\epsilon$ and
$v^\epsilon$ tend to $v$ weakly in $W^{1,2+\alpha}$ and then in
particular uniformyl on $[-R,R]$. We obtain that 
$$v=u.$$

This ends the proof.

 \begin{proposition}

 For $x\in]-R,R[$ let  

 $$G(x,u,u^\prime,u^{\prime\prime}):=a(x)|u'|^\alpha u'' + b(x)|u'|^\alpha u'$$
with $A\geq a(x)\geq a>0$, continuous on $[-R,R]$ and $b$ bounded,
  then there
exist some constants $C_1>0$, $C_2>0$ which depend on $a$ and the bound of $b$ such
that 

$$\bar\lambda> {C_1e^{-C_2R}\over R^{2+\alpha}}.$$

\end{proposition}
{\bf Proof}
Let $$B(x)=\int_{-R}^x \frac{b(x)(\alpha+1)}{a(x)}dx.$$
Then it is easy to show that 

 $$\bar\lambda\geq \lambda_1:=\inf_{u\in W_0^{1,2+\alpha}(]-R,R[)}\left\{ {\int_{-R}^R |u^\prime|^{\alpha+2}
e^{B(x)}dx\over \int_{-R}^R
{\alpha+1\over a(x)}|u|^{\alpha+2}e^{B(x)}dx}\right\}.$$
Indeed, the infimum is achieved and $u$, a function achieving the infimum,
 is a weak solution of 
$$|u^\prime|^\alpha (a(x)u^{\prime\prime}+b(x) u^\prime)= -\lambda_1
|u|^\alpha u.$$
Due to the previous lemma  $u$ is also a viscosity solution. 
One can assume that $u\geq 0$, so $u>0$ in $\Omega$,  using
strong maximum principle of Vazquez.   Hence, by
definition, $\bar\lambda\geq \lambda_1$.

But one has, for some universal constant $C$

$$\lambda_1 \geq { a\over \alpha+1}e^{-2|{b(x)\over a(x)}|_\infty R(\alpha+1)}
\inf_{u\in W_0^{1,2+\alpha} (]-R,R[)} {\int_{-R}^R|u^\prime|^{2+\alpha}\over
\int_{-R}^R |u|^{2+\alpha}}= {{C a\over \alpha+1}e^{-2|{b\over a}|_\infty R(\alpha+1)} \over
R^{2+\alpha}}$$
which is the desired result. 
This ends the proof.

\bigskip

\noindent {\em Proof of  Proposition \ref{propmin}.}

\noindent Suppose that $\Omega$ is contained in $[-R,R]\times \R^{N-1}$, let us define 
$$u(x) =3^qR^q-(x_1+2R)^q,$$where 
$q= {2.3^q R|b_1|_\infty \over a}+2$, then 
$|\partial_{x_1}u|\geq qR^{q-1}$ and 

$$a\partial_{x_1x_1}u+b_1 \partial_{x_1}u\leq -q{a(q-1)\over 2}R^{q-2}.$$
Finally, using also $u(x)\leq 3^q R^q,$ one gets

$$G(x, u, \nabla u, D^2 u)\leq -cq^{2+\alpha} R^{q(\alpha+1)-\alpha-2}\leq -c R  ^{-\alpha-2} 3^{-q(\alpha+1)} q^{2+\alpha}u^{1+\alpha},$$
and, by definition,

$$\bar\lambda\geq c R  ^{-\alpha-2} 3^{-q(\alpha+1)} q^{2+\alpha}.$$
Using the expression of $q$ in function of $b_1$ one gets the announced   estimate. 

\begin{remark}
Let us note that in the case $b = cte $ or when there exists  some
direction $e_1 $ such that $b(x).e_1= cte$ and $\Omega$ is bounded in
this direction one has a better estimate.

Indeed, similarly to \cite{BD2} one can consider 
$$u(x) = u(x_1) = 7R^2-x_1^2-3(sign b_1)Rx_1.$$
This function is positive on $x_1\in ]-R,R]$,
its  gradient is never zero. Hence one has, for some constant $C$: 

\begin{eqnarray*}
G(x,u,\nabla u, D^2 u) &= &|3(sign b_1)R +2x_1|^\alpha (-2+ b_1(-3R(\rm{sign}
b_1)-2x_1))\\
&\leq& C R^\alpha (-2-3|b_1|R+2|b_1|R)\\
&\leq& C R^\alpha (-2-|b_1|R).
\end{eqnarray*}
This implies that $$G(x,u,\nabla u, D^2u)\leq -C(2+|b_1|R) R^{-2-\alpha}u^{1+\alpha}$$
which yields
$$\bar \lambda\geq {C_1\over R^{2+\alpha}} + {C_2|b_1|\over
R^{1+\alpha}}$$ which is a more accurate  lower bound than in the general
case.
\end{remark}

\begin{proposition}\label{pro1}

Suppose that $R$ is the radius of the largest ball contained in  
$\Omega$ and suppose that $F$ satisfies assumption  (H1) and (H2).

Furthermore let $b$ and $c$ be bounded functions.  Then, there exists
some constant $C_1$ which depends only on
$N$, $\Omega$ $\alpha$, $a$ and $A$,  such that 
$$\bar\lambda\leq C_1\left({1\over R^{\alpha+2}}+{|b|_\infty\over 
R^{\alpha+1}}\right)+| c|_\infty.$$

\end{proposition}
{\bf Proof.} Without loss of generality we can suppose that the largest ball contained in $\Omega$
is $B_R(0)$. 
Let $\sigma$ be defined as 
$$\sigma(x) = {1\over 2q} (|x|^q-R^q)^2$$
with $q = {\alpha+2\over \alpha+1}$, for $x\in B_R(0)$. 

We need to compute the supremum in $B_R(0)$ of 
$$\frac{
-F(x,\nabla \sigma, D^2\sigma) -b(x).\nabla \sigma |\nabla
\sigma|^\alpha}{\sigma^{\alpha+1}}.$$

Let $\sigma(x)=g(r) $, for $r=|x|$. Clearly
$g^\prime (r) = r^{2q-1}-r^{q-1} R^q$ and
$$g^{\prime\prime }(r) = (2q-1) r^{2q-2}-(q-1) r^{q-2} R^q.$$

Furthermore $g'\leq 0$ while $g''\leq 0$ for $r\leq
\left(\frac{q-1}{2q-1}\right)^{1\over q}R$ and positive
elsewhere. Hence by condition (H2) and using the fact that for radial
functions the eigenvalues of the Hessian are
$\displaystyle{\frac{g'}{r}}$ with multiplicity N-1 and $g''$ (see \cite{CuL}) we get

-for $r\leq {\left(\frac{q-1}{2q-1}\right)^{1\over q}R}$   

$$ |F(x,\grad \sigma,D^2\sigma)|
\leq |g'|^\alpha  \left\vert ag^{\prime\prime}(r)+ a({N-1\over r})
g^\prime(r)\right\vert
$$

while for $r\geq \left(\frac{q-1}{2q-1}\right)^{1\over q}R$

$$- F(x,\grad \sigma,D^2\sigma) \leq
-|g'|^\alpha\left[a g^{\prime\prime}(r)+
A({N-1\over r}) g^\prime(r)\right] 
$$

i.e.

$$
-F(x, \grad \sigma,D^2\sigma) \leq |g'|^\alpha r^{q-2} (-B_1r^q+B_2R^q)
$$
where  $B_1=a(2q-1)+A(N-1)$ and $B_2=a(q-1)+A(N-1)$.

Let $R_1$ be defined as 
$$R_1 = R\left({B_2+ |b|_\infty R^{(q-1)(\alpha+1)}\over
 B_1+ |b|_\infty R^{(q-1)(\alpha+1)}}\right)^{\frac{1}{q}}< R$$ 
then for $r\geq R_1$
$$-F(x,\nabla \sigma, D^2\sigma) -b(x).\nabla \sigma |\nabla
\sigma|^\alpha \leq 0.$$
Hence the supremum is achieved for $r\leq R_1$. On that set one can
use  an upper bound for 
$|-F(x,\grad \sigma,D^2\sigma)+ b.\nabla \sigma |\nabla \sigma|^\alpha |$
and a lower bound for $\sigma$ e.g.
$$\sigma \geq \frac{1}{2q}(R^q-R_1^q)^2.$$

More precisely for $r\leq 
R_1$ and for some constants $C_1$, $C_2$, $C_1^\prime$ $C_2^\prime$ 
depending on
$a, A, N$ and $|b|_\infty$
one has:
\begin{eqnarray*}
\frac{|-F(x,\grad \sigma,D^2\sigma)- b.\nabla \sigma |\nabla \sigma|^\alpha |}{\sigma^{\alpha+1}} &
\leq &  r^{q(\alpha+1)-\alpha-2} {C_1R^q\over (R^q-R_1^q)^{\alpha+2}}+ {C_2|b|_\infty
R^{(q-1)(\alpha+1)}\over R^{q(\alpha+1)}}\\
&\leq& {C^\prime _1\over
R^{\alpha+2}}+ {C_2^\prime\over R^{\alpha+1}}
\end{eqnarray*}

Then $\sigma$ is a subsolution in $B_R(0)$ of

$$F(x,\nabla\sigma,D^2 \sigma ) +b.\nabla \sigma |\nabla
\sigma|^\alpha +c|\sigma|^{\alpha}\sigma+ ({C_1\over R^{\alpha+2}}+
{C_2\over R^{\alpha+1}}+|c|_\infty )|\sigma| ^{\alpha}\sigma 
 \geq 0,$$
with $\sigma=0$ on $\partial B_R(0)$.
Suppose by contradiction that $\bar\lambda > {C_1\over R^{\alpha+2}}+{C_2\over R^{\alpha+1}} 
+|c|_\infty$.  Clearly since $B_R(0)\subset\Omega$, $\bar\lambda(B_R(0))\geq\bar\lambda> {C_1\over R^{\alpha+2}}+{C_2\over R^{\alpha+1}} 
+|c|_\infty$, and then according to the maximum principle, Theorem \ref{maxp}, one should have  that 
$\sigma \leq 0$ in $B_R(0)$, a contradiction.
This ends the proof.

\subsection {Comparison theorem for $\lambda<\bar\lambda$}

\begin{theorem}\label{complambda}
Suppose that $F$ satisfies (H1), (H2'), and (H4),  
that $b$ and $c$ are continuous and
bounded and $b$ satisfies
$(H5)$.  Suppose that $\tau< \bar\lambda$, $f\leq 0$, $f$ is
upper semi-continuous and
$g$ is lower semi-continuous  with $f\leq  g$.

Suppose that there exist
$\sigma$ upper semi continuous , and
$v$ non-negative and lower semi continuous , satisfying
\begin{eqnarray*}
F(x, \grad v,D^2v)+b(x).\nabla v |\nabla
v|^{\alpha}+(c(x)+\tau) v^{1+\alpha} & \leq & f \quad  \mbox{in}\quad \Omega \\
F(x,  \grad \sigma,D^2\sigma) +b(x).\nabla
\sigma |\nabla \sigma |^{\alpha} +(c(x)+\tau)|\sigma|^{\alpha}\sigma
& \geq & g \quad  \mbox{in}\quad \Omega  \\ 
\sigma \leq  v &&   \quad  \mbox{on}\quad \partial\Omega 
\end{eqnarray*}
Then $\sigma\leq v$ in $\Omega$ in each of these two cases:

\noindent 1) If $v>0$ on $\overline{ \Omega}$ and either $f<0$ in 
$\Omega$, 
 or  $g(\bar x)>0$ on every point $\bar x$ such that  $f(\bar x)=0$,

\noindent  2) If $v>0$ in $\Omega$, $f<0$ and $f<g $ on
$\overline\Omega$ 
\end{theorem}

{\bf Proof.}
We act as in the proof of Theorem 3.6 in \cite{BD1}.

1) We assume first that $v>0$ on $\overline{\Omega}$. 
Suppose by contradiction that
$\sigma > v$ somewhere in
$\Omega$. The supremum of the function $\displaystyle{\sigma\over v}$  on
$\partial
\Omega$ is  less  than $1$ since $\sigma\leq v$ on $\partial \Omega$ and
$v>0$ on
$\partial \Omega$, then its supremum is achieved inside
$\Omega$. Let
$\bar x$ be a point such that 
$$1<{\sigma(\bar x)\over v(\bar x)}= \sup_{x\in \overline{\Omega}} {\sigma(x)\over v(x)}.$$

We define 
$$\psi_j(x,y) = {\sigma(x)\over v(y)}-{j\over qv(y)} |x-y|^q.$$
For $j$ large enough, this function achieves its maximum which is greater than 1, on
some couple
$(x_j, y_j)\in
\Omega^2$. It is easy to see that this sequence converges to
$(\bar x,
\bar x)$,  a maximum point for ${\sigma\over v}$.

\medskip
Since on test functions that have zero gradient the definition of viscosity solutions doesn't require to test equation,  we need to prove first 
 that  $x_j$, $ y_j$ can be chosen such that $x_j\neq y_j$ for $j$ large enough.   

Indeed, if $x_j=y_j$ one would have 
for all $x\in \Omega$ 
$${\sigma (x)-{j\over q} |x_j-x|^q\over v(x_j)}\leq C = {\sigma (x_j)\over
v(x_j)},$$ which implies that 
$$\sigma (x)\leq \sigma (x_j)+{j\over q} |x_j-x|^q.$$
This means that $\sigma$ has a local maximum on the point $x_j$. We
argue as it is done in the proof of Theorem \ref{compprnew} :
If $x_j$ is not a strict local maximum then $x_j$ can be replaced by
$x^\prime_j$ close to it and then $(x^\prime_j, x_j)$ is also a maximum
point for $\psi_j$.

 If the maximum is strict using  
Lemma \ref{dem22} one gets that 
$$(c(x_j) +\tau)\sigma(x_j)^{1+\alpha} \geq g(x_j).$$
 But one also has  
$$v(x)\geq v(x_j)-{j\over q C} |x_j-x|^q.$$
hence $x_j$ is a local minimum for $v$, and if it is not strict, there
exists $x_j^\prime$ which is different from $x_j$ such that $(x^\prime_j,
x_j)$ is also a maximum point for $\psi_j$. 

If the minimum is strict  
using once more Lemma \ref{dem22} one would have 
$$(c(x_j)+\tau) v^{1+\alpha}(x_j) \leq f(x_j).$$

This is a contradiction for $j$ large enough. Indeed,  passing
to the limit one would get 
$$(c(\bar x)+\tau) (\sigma (\bar x)^{1+\alpha}-v(\bar x)^{1+\alpha}) \geq
g(\bar x)-f(\bar x)\geq 0.$$
Since $\sigma(\bar x)>v(\bar x)$ this implies that 
$$c(\bar x)+\tau\geq 0.$$
Now there are two cases either  $f(\bar x)<0$ or $f(\bar x)=0$  and the above inequality is strict.
In both cases it contradicts

$$(c(\bar x)+\tau) v(\bar x)^{1+\alpha} \leq f(\bar x).$$

We can take $x_j$ and $y_j$ such that $x_j\neq y_j$. 

\bigskip
Moreover there exist
$X_j$ and $Y_j$ such that 

$$\left(j|x_j-y_j|^{q-2} (x_j-y_j), {X_j\over v(y_j)}\right)\in J^{2,+} \sigma(x_j)$$
and 

$$\left(j|x_j-y_j|^{q-2} (x_j-y_j){v(y_j)\over \beta_j}, {-Y_j\over \beta_j}\right)\in J^{2,-} v(y_j)$$
 where $\beta_j = \sigma(x_j)-{j\over q} |x_j-y_j|^q$
and 
$$F(x_j, j|x_j-y_j|^{q-2} (x_j-y_j), X_j)- F(y_j, j|x_j-y_j|^{q-2}
(x_j-y_j), {-Y_j})\leq \omega (v(y_j)j|x_j-y_j|^q).$$

We can use the fact that $\sigma$ and $v$ are respectively sub and super solution to obtain:

\begin{eqnarray*}
g(x_j)-\tau \sigma (x_j)^{1+\alpha} &-&c(x_j) \sigma
(x_j)^{1+\alpha}\leq F(x_j, j|x_j-y_j|^{q-2} (x_j-y_j), {X_j\over
v(y_j)})\\
&&+b(x_j).
j^{1+\alpha}|x_j-y_j|^{(q-1)(1+\alpha)-1} (x_j-y_j)\\
 &\leq & {\beta_j^{1+\alpha}\over v(y_j)^{1+\alpha}} 
\Bigl\{F(y_j, j|x_j-y_j|^{q-2}  (x_j-y_j){v(y_j)\over \beta_j}, {-Y_j\over
\beta_j})\\
&&+\omega(jv(y_j)|x_j-y_j|^q)\Bigr\}+\\
&& + b(x_j).
j^{1+\alpha}|x_j-y_j|^{(q-1)(1+\alpha)-1} (x_j-y_j)\\ 
&\leq&
{\beta_j^{1+\alpha}\over v(y_j)^{1+\alpha}} 
\Bigl\{F(y_j, j|x_j-y_j|^{q-2}  (x_j-y_j){v(y_j)\over \beta_j}, {-Y_j\over
\beta_j})\\
&&+b(y_j).
j^{1+\alpha}|x_j-y_j|^{(q-1)(1+\alpha)-1}
(x_j-y_j)\left({v(y_j)\over
\beta_j}\right)^{1+\alpha}\Bigr\}\\
&& +{\omega(v(y_j)j|x_j-y_j|^q)\over
v(y_j)^{1+\alpha}} +C(j|x_j-y_j|^q)^{1+\alpha}\\
 &\leq &
(-\tau-c(y_j)) \beta_j^{1+\alpha} + {\beta_j^{1+\alpha}\over
v(y_j)^{1+\alpha}} f(y_j)+o(1).
\end{eqnarray*}
Passing to the limit, since $c$ is continuous, we get:
$$g(\bar x)\leq  \left(\frac{\sigma(\bar x)}{v(\bar x)}\right)^{\alpha+1} f(\bar x).$$
Either $f(\bar x)=0$ and then we have reached a contradiction because, in that case, by hypothesis
$$g(\bar x)>0,$$
or $f(\bar x)<0$,
and then we get
$$0<f(\bar x)\left[1-\left(\frac{\sigma(\bar x)}{v(\bar x)}\right)^{\alpha+1}\right]\leq f(\bar x)-g(\bar x)\leq 0.$$
This concludes the proof of the first part.

\bigskip

2) For the second part, let $m$ be  such that $f-g\leq -m<0$, and
$f< -{m\over 2}$. Let
$\epsilon$ be given such that by the uniform continuity of the function 
$(x+\epsilon)^{1+\alpha}$ on $[0, |v|_\infty]$ one has 
$$|\lambda +c|_\infty\cdot|(v+\epsilon)^{1+\alpha}-v^{1+\alpha}|\leq {m\over 2}.$$
Then $w=v+\epsilon$ is a  supersolution of 
$$F(x,\grad w, D^2w)+(\lambda+c) w^{1+\alpha} \leq f+{m\over 2}
\leq g-{m\over 2}<g
\leq F(\sigma,\grad \sigma,D^2\sigma) + (\lambda+c)
|\sigma|^\alpha\sigma.$$ We are now in a position to use the first part
of the theorem, since 
$$w=v+\epsilon>0\quad\mbox{and}\quad u\leq v\leq v+\epsilon \quad\mbox{on}\quad \partial\Omega.$$
 and then 
$u\leq v+\epsilon$
in $\Omega$. Letting $\epsilon$ go to zero we get the required conclusion.
This ends the proof.

\section{Regularity results}

In this section we shall prove that the viscosity solutions are H\"older continuous. Since the H\"older estimates depend only on the bounds of $f$ and the structural constants, this H\"older continuity will allow us to have a compactness criteria that will be useful in the next section.
Let us note that we state all the results with $c=0$. Indeed, one can
consider $c(x)|u|^\alpha u$ in the right hand side since it is bounded, and get
the same regularity results. 

\begin{proposition}
Suppose that $F$ satisfies (H1), (H2), (H3). 
Let $f$ be a bounded function in $\overline{\Omega}$.
 Let $u$ be a viscosity non-negative bounded solution of 

\begin{equation}\label{eq4.1}
\left\{
\begin{array}{lc}
F(x, \nabla u, D^2u)+ b(x).\nabla u|\nabla u|^{\alpha}=f & \ {\rm in}\ \Omega\\
u=0 &  \ {\rm on}\ \partial\Omega.
\end{array}
\right.
\end{equation}
Then if $b$ is bounded,  for
any
$\gamma\in (0,1)$ there exists some constant
$C$ which depends only on $|f|_\infty$ and $|b|_\infty$ such that for any $(x,y)\in\bar\Omega^2$

$$|u(x)-u(y)|\leq C|x-y|^\gamma.$$
\end{proposition}
An immediate consequence of the above Proposition is the 
\begin{corollary} \label{comp}
Suppose that $F$ satisfies (H1), (H2) and (H3). Suppose that $f_n$ is a  
sequence of continuous and uniformly bounded functions, and $u_n$ is a sequence of bounded viscosity solutions of  

$$F(x, \nabla u_n, D^2u_n)+ b(x).\nabla u_n|\nabla
u_n|^{\alpha}=f_n(x)$$ 
with $b$ bounded, $u_n=0$ on $\partial\Omega$. Then the sequence  $u_n$ is relatively
compact in
${\mathcal C} (\overline{\Omega})$.
\end{corollary}
{\bf Proof.}
The proof  relies on ideas
used to prove H\"older and Lipschitz  estimates in \cite{IL}, as it is
done in \cite{BD2}.

We use Proposition \ref{ah} in section 3 which implies in particular
that there exists $M_0$ such that  
\begin{equation}\label{claim}
u(x) \leq M_0d(x)^\gamma
\end{equation}
for $d(x):=d(x,\partial\Omega)\leq \delta$.

\bigskip
We now prove H\"older's regularity inside $\Omega$. 

We construct a function $\Phi$ as follows:
Let $M_o$ and $\gamma$ be as in (\ref{claim}),
  $M= \sup  (M_o, {2\sup u\over \delta^\gamma})$ and 
 $\Phi (x) = M |x|^\gamma$. We also define  
$$\Delta_\delta=\{(x,y)\in \Omega^2,\ |x-y|<\delta\}.$$

\noindent {\bf Claim} {\em For any} $(x,y)\in\Delta_\delta$
\begin{equation}\label{ho}
u^\star(x)-u_\star(y)\leq \Phi(x-y).
\end{equation}

If the Claim  holds this completes the proof, indeed
taking $x=y$ we would get that $u^\star=u_\star$ and then $u$ is
continuous.  Therefore, going back to (\ref{ho}),
$$u(x)-u(y)\leq {2\sup u\over \delta^\gamma} |x-y|^\gamma,$$ for
$(x,y)\in \Delta_\delta$
 which is equivalent to the local H\"older continuity. 

\bigskip

\noindent Let us check  that  (\ref{ho}) holds 
 on $\partial \Delta_\delta$. On that set:

- either $|x-y|=\delta$ and then
$u^\star(x)-u_\star(y)\leq M\delta^\gamma$ 
since  $M\delta^\gamma \geq {2\sup u}$,  

-or
$(x,y)\in\partial(\Omega \times \Omega)$. In that case, for 
$(x,y)\in (\Omega\times\partial \Omega)$ we have just proved that
$$u^\star(x)\leq M_od(x)^\gamma\leq M|y-x|^\gamma,$$
while for $(x,y)\in \partial \Omega\times \Omega $
$$0-u_\star (y)\leq 0\leq M_0 |y-x|^\gamma.$$

Now we consider interior points. Suppose by
contradiction that $u^\star(x)-u_\star(y)> \Phi(x-y)$
for some $(x,y)\in \Delta_\delta$. Then there exists $(\bar x, \bar y)$ such that 
$$u^\star(\bar x)-u_\star(\bar y)-\Phi(\bar x-\bar y)=\sup (u^\star(x)-u_\star(y)-\Phi(x-y))>0.$$
 Clearly $\bar x\neq \bar y$.
Then using Ishii's Lemma \cite{I},  there exist $X$ and $Y$ such that
$$( \gamma M(\bar x-\bar y)|\bar x-\bar y|^{\gamma-2}
, X)\in J^{2,+ }u^\star(\bar x)$$
$$(\gamma M(\bar x-\bar y)|\bar x-\bar y|^{\gamma-2}
, -Y)\in J^{2,-}u_\star(\bar y)$$
with 
$$\left(\begin{array}{cc}
X&0\\
0&Y
\end{array}\right) \leq \left(\begin{array}{cc}
B &-B\\
-B&B
\end{array}\right) $$
and  $B = D^2\Phi(\bar x-\bar y)$.

 We need a 
more precise estimate, as in
\cite{IL}. For that aim let :

$$0\leq P : = {(\bar x-\bar y\otimes \bar x-\bar y)\over |\bar x-\bar y|^2}\leq I.$$

Using $-(X+Y)\geq 0$ and $(I-P)\geq 0$ and the properties of the
symmetric matrices one has 
$$tr(X+Y)\leq tr(P(X+Y)).$$ Remarking in addition that
$X+Y\leq 4B$, one  sees that
$tr(X+Y)\leq tr(P(X+Y))\leq 4tr(PB)$. But $tr (PB)=\gamma M(\gamma-1)|
\bar x-\bar y|^{\gamma-2}<0$, hence

\begin{equation}\label{elip}
|tr (X+Y)|\geq 4\gamma M(1-\gamma)| \bar x-\bar y|^{\gamma-2}.
\end{equation}

Furthermore  by Lemma III.1 of \cite{IL} there exists a universal constant $C$ such that 
$$|X|, |Y|\leq C (|tr(X+Y)|+ |B|^{1\over 2} |tr(X+Y)|^{1\over
2})\leq C|tr(X+Y)|,$$
since $|B|$ and $|tr(X+Y)|$ are of the same order. 
Now we can use the fact that $u$ is both a sub and a super solution of 
(\ref{eq4.1}), 
and applying condition (H2), (H3) concerning $F$ :

\begin{eqnarray*}
f(\bar x)&-&  (\gamma M)^{1+\alpha} b(\bar x).(\bar x-\bar y)|\bar x-\bar
y| ^{(\gamma-1)(\alpha+1)-1}\\
&\leq&   F(\bar x,\grad\Phi(\bar x-\bar y),X)\\
&\leq & F(\bar y, \grad\Phi(\bar x-\bar y), X)+\tilde \omega(|\bar x-\bar
y|) |\nabla \Phi(\bar x-\bar y) |^\alpha|X|\\
 &\leq &  a|\grad\Phi(\bar x-\bar y) |^\alpha tr (X+Y) + 
F(\bar y, \grad_y\Phi(\bar y), tr(-Y))\\
&&+\tilde \omega(|\bar x-\bar y|)
|\nabla \Phi(\bar x-\bar y) |^\alpha |X|\\ &\leq &
f(\bar y)    + a|\grad\Phi(\bar x-\bar y) |^\alpha tr (X+Y)+ (\gamma
M)^{1+\alpha} |b|_\infty |\bar x-\bar y|^{(\gamma-1)(\alpha+1)}\\
&&+ \tilde \omega(|\bar x-\bar y|)
|\nabla \Phi(\bar x-\bar y) |^\alpha |tr(X+Y)|.
\end{eqnarray*}
Which implies, using (\ref{elip}) 
$$ |\grad\Phi(\bar x-\bar y) |^\alpha\gamma M(1-\gamma)| \bar x-\bar
y|^{\gamma-2}\left(a-C\tilde \omega(|\bar x-\bar y|)-2{|b|_\infty\over(1-\gamma)} |\bar x-\bar
y|\right)\leq f(\bar y)-f(\bar x).
$$ We choose $\delta$ small enough in order that 
$C\tilde \omega (\delta)+2{|b|_\infty\over(1-\gamma)} \delta <{a\over 2}$. 
Recalling that
$|\grad\Phi(\bar x-\bar y) |=\gamma M|\bar x-\bar y|^{\gamma-1}$ the
previous inequality becomes: 
\begin{equation}\label{M}
 {a\over 2}M^{\alpha+1}\gamma^{1+\alpha}(1-\gamma)|\bar x-\bar y|^{\gamma
(\alpha+1)-(\alpha+2)}\leq 2 |f|_\infty.
\end{equation} 
Using $M\geq {2(\sup u)\over \delta ^\gamma} $ and
$|\bar x-\bar y|\leq \delta$ one obtains 

$$a(2\sup u)^{1+\alpha} \gamma ^{1+\alpha} (1-\gamma)
\delta^{-(\alpha+2)}
\leq 4|f|_\infty. $$ This is clearly false for $\delta$
small enough and it concludes the proof.
This ends the proof.

\bigskip
For completeness sake we shall now prove some Lipschitz  regularity of the solution.
To get Lipschitz regularity we need a further assumption as it was done in \cite{BD2}. 
Let us remark that Lipschitz regularity is not necessary to prove 
the existence results, hence this further assumption will be used only in the present\ part of the
paper.

(H7){\it
There exists $\nu>0$ and $\kappa \in ]1/2,1]$ such that for all $|p|=1$ ,
$|q|\leq {1\over 2}$, $B\in {\mathcal S}$
 $$|F(x,p+q,B)- F(x,p,B)|\leq \nu |q|^\kappa |B|$$
which implies by homogeneity 
that for all } $p\neq 0$ ,
$|q|\leq {|p|\over 2}$, $B\in {\mathcal S}$
 $$|F(x,p+q,B)-F(x,p,B)|\leq \nu |q|^\kappa|p|^{\alpha-\kappa} |B|$$

One has, then, the following regularity result: 
\begin{theorem}\label{L}
If $F$ satisfies (H1),(H2), (H3) and (H7) and  if $b$ is bounded, then the
bounded solutions of
$$\left\{\begin{array}{lc}
F(x,\nabla u,D^2 u)+b(x).\nabla u|\nabla u|^\alpha = f & {\rm in}\ \Omega\\
u=0 & {\rm on}\ \partial \Omega
\end{array}
\right.
$$
are Lipschitz continuous inside $\Omega$. 
\end{theorem}

\noindent {\it Proof of Theorem \ref{L}.}  The proof proceeds similarly to
the proof given by Ishii and Lions in \cite{IL} and as it is required
in that paper, we use the fact that we already know that $u$ is H\"older
continuous, together with the additional assumption (H7) .

To simplify the calculation but, without loss of generality we shall suppose that in hypothesis (H2)
  $a= A =1$.  Let $\gamma$ be in $]{1\over 2\kappa},1[$ and $c_\gamma$
such that by the H\"older's continuity  proved before
$$|u(x)-u(y)|\leq c_\gamma |x-y|^\gamma.$$

Let $\mu$ be an increasing function such that $\mu(0)=0$, and
$\mu(r)\geq r$, let $l(r) = \int_0^r ds\int_0^s {\mu(\sigma)\over
\sigma} d\sigma$, let us note that since $\mu\geq 0$, for $r>0$:
$$l(r)\leq rl^\prime (r)$$

Let $r_0$ be such that $l^\prime (r_0) = {1\over 2}$. Let also $\delta>0 $ be given, $K= {r_0\over
\delta}$, and $z$ be such that $d(z, \partial \Omega)\geq 2\delta$. 

We define  $\varphi(x,y) =
\Phi(x-y)+L|x-z|^k$
where $\Phi(x) =M(K|x| -l(K|x|))$, 
and  
$$\Delta_z= \{ (x,y)\in \R^N\times \R^N, |x-y|<\delta, |x-z|< \delta\}.$$
We shall now choose all the constants above.

- $k$ is such that   $k={1\over 1-{1\over 2\gamma\kappa }}$

 - $M$ and $L$ are such that  $M= {4\sup u\over r_0}$ and 
$L= c_\gamma\delta^{\gamma-k}$, using the H\"older continuity of $u$,
 one has  
$$u(x)-u(y)\leq \varphi(x,y)$$ on $\partial \Delta_z$. 
Indeed, the assumption on $r_0$ implies that $\Phi(x)\geq MK {|x|\over
2}$ for $|x|\leq r_0$ and then if $|x-y|=\delta$, 
\begin{eqnarray*}
u(x)-u(y)&\leq& 2\sup u
\leq  {Mr_0\over 2}\leq {MK\delta\over 2}\\
&\leq &\Phi(x-y)
\leq \varphi(x,y),
\end{eqnarray*}
while if $|x-z|=\delta$
$$u(x)-u(y)\leq c_\gamma |x-y|^\gamma \leq c_\gamma \delta^\gamma =
L|x-z|^k \leq \varphi(x,y).$$

Suppose by contradiction that for some point $(\bar x, \bar y)$ one has 
$$u(\bar x)-u(\bar y)> \varphi(\bar x, \bar y).$$
Clearly $\bar x\neq \bar y$. 
Note that 
$$L|\bar x-z|^k\leq c|\bar x-\bar y|^\gamma.$$
Proceeding as in the previous proof, there exist 
$X$, $Y$ such that 

$$( MK(\bar x-\bar y)|\bar x-\bar y|^{-1}
(1-l^\prime (K|\bar x-\bar
y|))+ kL |\bar x-z|^{k-2} (\bar x-z), X)\in J^{2,+}u(\bar x),$$
and

$$( MK{\bar x-\bar y\over |\bar x-\bar y|}(1-l^\prime (K|\bar x-\bar
y|)), -Y)\in J^{2,-}u(\bar y),$$
where the matrices $X$ and $Y$ satisfy
\begin{equation}\label{ish}
\left(\begin{array}{cc}
X&0\\
0&Y
\end{array}\right) \leq \left(\begin{array}{cc}
B+\tilde L&-B\\
-B&B
\end{array}\right) 
\end{equation}
with $B = D^2\Phi(\bar x-\bar y)$
and

 $$\tilde L = kL|\bar x-z|^{k-2} \left(I+ (k-2){(\bar x-z\otimes \bar x-z)\over
|\bar x-z|^2}\right).$$
Let us note that similarly to the H\"older case, (\ref{ish}) implies that
$X+Y-\tilde L\leq 4B$ and then 

$$tr(X+Y-\tilde L )\leq 4 tr(PB)$$
with 

$$P = {((\bar x-\bar y)\otimes (\bar x-\bar y))\over |\bar x-\bar y|^2}.$$
This gives: 

$$tr(X+Y-\tilde L)\leq -{MK \mu(K|\bar x-\bar y|))\over |\bar x-\bar
y|}\leq -MK^2.$$

 Let us note that  
$$\nabla_x \varphi(x) = MK(1-l^\prime (K|\bar x-\bar y|)){\bar x-\bar y\over
|\bar x-\bar y|} +k L|\bar x-z|^{k-2} (\bar x-z) $$
and 

$$L|\bar x-z|^{k-1} =O(\delta^{\gamma-k}\delta^{k-1})=O(K^{1-\gamma}).$$
From this we get in particular that for $\delta>0$ small enough (or $K$
large enough)
$${2MK}\geq (|\nabla_x \varphi(\bar x, \bar y)|, |\nabla_y \varphi(\bar x,
\bar y)| ) \geq {MK\over 4}.$$ Finally observe that 
$|\tilde L|\leq k(k-1)L |\bar x-z|^{k-2} \leq (C\delta^{\gamma-k})^{2\over
k}  (\delta)^{\gamma (k-2)/k} = O (\delta ^{\gamma-2}) = O(K^{2-\gamma})$,
from which we derive  that for $K$ large enough 
$tr(X+Y) \leq 0$ and 
$$|tr(X+Y)|\geq C(K^2)$$
for some $>0$ universal constant $C$, and 
$|\tilde L|\leq |tr(X+Y)|$
for $K$ large enough.  

In the following we shall need a bound from above for  $|X|$.  In
order to make the reading easier the constants $C$ or $c$ will be constants
which depend only on the data, and they may vary from one line to
another. 
Remark that the lemma III.1 in \cite {IL} ensures the existence of some
universal constant such that 
$$|X-\tilde L|+|Y|\leq C\left(|B|^{1\over 2} |tr(X+Y-\tilde L)|^{1\over 2}
+|tr(X+Y-\tilde L)|\right)$$
with 
$B= D^2\varphi$,  and with the considerations on $\tilde L
$ with respect to $|tr(X+Y)|$
one also has
$$|X|+|Y|\leq C\left(|B|^{1\over 2} |tr(X+Y)|^{1\over 2}
+|tr(X+Y)|\right).$$

Let us note that 
$$|D^2\varphi|\leq {CK\over |\bar x-\bar y|},$$
and then  with the assumptions on $\mu$,  
$|tr(X+Y)|\geq C\geq K|tr(X+Y)|^{1\over
2}$ from which one derives that   
$$|X|\leq |tr(X+Y)|(1+{1\over K^{1\over 2} |\bar x-\bar y|^{1\over 2}}).$$
We need to prove that 
$$|\nabla_x \varphi(\bar x, \bar y)|^{\alpha-\kappa}(\tilde L|\bar
x-z|^{k-1} )^{\kappa} |X|= o (|tr(X+Y)||\nabla \varphi|^\alpha).$$
For that aim we write 
\begin{eqnarray*}
|\nabla _x\varphi(\bar x, \bar y)|^{\alpha-\kappa} ( L|\bar
x-z|^{k-1} )^{\kappa} |X|&\leq &
cK^{\alpha-{\gamma\kappa\over
k}}|\bar x-\bar y|^{\gamma\kappa (1-{1\over k})}|X|\\
&\leq & 
 cK^{\alpha-{\gamma\kappa\over k} }|\bar x-\bar
y|^{1\over 2} |tr(X+Y)|(1+{1\over K^{1\over 2} |\bar x-\bar y|^{1\over
2}})\\
&\leq &c|tr(X+Y)| (K^{\alpha-{\gamma\kappa\over k}-{1\over 2}})\\
&= &c|tr(X+Y)| K^{\alpha} K^{-\gamma\kappa-{1\over 2}}= o(|tr(X+Y)
|\nabla
\varphi|^\alpha)
\end{eqnarray*}

We now obtain using assumption (H2) and (H3 ) concerning $F$
\begin{eqnarray*}
f(\bar x)-b(\bar x).\nabla_x\varphi|\nabla_x \varphi|^\alpha &\leq&
F(\bar x, \nabla_x
\varphi(\bar x,
\bar y), X)+|b|_\infty O(K^{1+\alpha})\\
 &\leq & F(\bar y,\nabla_y\varphi(\bar x, \bar y), X)+\nu |L|\bar
x-z|^{k-1}|^{\kappa} |\nabla _x\varphi|^{\alpha-\kappa}| |X|\\
& &+|b|_\infty
O(K^{1+\alpha})\\ &\leq& F(\bar y, \nabla_y\varphi(\bar x, \bar y),
-Y)+|\nabla
\phi|^\alpha |tr(X+Y) |\\
& &+b(\bar y).\nabla_y
\varphi|\nabla_y\varphi|^{\alpha}+O(K^{-{\gamma\kappa\over k}-{1\over 2}}
|\nabla
\varphi|^\alpha)  tr(X+Y)\\
& &+|\nabla
\phi|^\alpha  tr(X+Y)+2|b|_\infty O(K^{1+\alpha})
\\ 
&\leq & f(\bar y) +O(K^{2+\alpha-\gamma\kappa})-C
(K^{\alpha+2})+|b|_\infty O(K^{1+\alpha})
\end{eqnarray*}
From this one gets a contradiction for $K$ large. 
We have proved that for all $x$ such that 
$d(x,\partial \Omega)\geq 2\delta$ and for $y$ such that $|x-y|\leq \delta$ 
$$u(x)-u(y)\leq \left({2\sup u\over r_0}\right)\left( {|x-y|\over
\delta}\right).$$ Recovering the compact set $\Omega$ by a finite number
of 
${\mathcal C}^2$ sets
$\Omega_i$, $\Omega_i\subset \Omega_{i+1}$  such that
$d(\partial \Omega_i, \partial \Omega_{i+1})\leq 2\delta$, the local
Lipschitz continuity is proved. 

\section{Existence's results}

We now prove the existence of non negative solutions of 

$$\left\{\begin{array}{lc}
G(x, u,\nabla u, D^2
u) +\lambda u^{1+\alpha} = -f &\ {\rm in}\ \Omega\\
u=0 &\ {\rm on}\ \partial\Omega
\end{array}
\right.
$$ where $f$ is a given positive function. The steps are the following :

Step 1: Exhibit a sub and a super solution
  of the equation when the coefficient of the zero order  is non positive  and $f$ is 
constant.

Step 2:  Under the same conditions on the zero order term, use Perron's method to solve the equation for any negative function
$-f$.

Step 3: From the previous steps
 we construct a solution of the above Dirichlet problem when
$\lambda<\bar\lambda$ without conditions on the sign of $c(x)$.

Step 4: This will also allow to prove the existence of the associated eigenvalue.

The   first step is obtained by remarking that $0$ is a sub solution and
establishing the following
\begin{proposition}\label{propbeta}
Suppose that $F$ satisfies (H1) and (H2), $b$  and $c$ are  bounded; furthermore let $c$ be 
non-positive in
$\overline\Omega$. Then there exists a function $u$ which
is a nonnegative viscosity super solution of 
$$\left\{\begin{array}{cc}
F(x, \grad u,D^2u)+ b(x) |\nabla u|^\alpha \nabla u+ c(x) u^{1+\alpha}\leq
-1&\ {\rm in}
\
\Omega\\ u=0&\ {\rm on}\ \partial \Omega
\end{array}
\right.
$$

\end{proposition}
{\bf Proof.}
Let $d$ be the distance function to $\partial \Omega$, which is well defined in
$\Omega$ and satisfies the properties stated in Remark \ref{dist}. 
Let  $K> diam \Omega$. Then $d\leq K$. Let $\gamma \in ]0,1[$ and let
$k$ be a large enough constant to be chosen later.
Let $u$ be defined as 

$$u(x) = 1-{1\over (1+d(x)^\gamma)^k}.$$
Clearly 
$u=0$ on the boundary. 

Suppose that $\psi$ is a ${\mathcal C}^2$ function  such that 
$(u-\psi)(x)\geq (u-\psi) (\bar x)=0$, for all $x$ in a small
neighbourhood of $\bar x$.  Then
$J^{2,-} u(\bar x)\neq
\emptyset$ and
 then the function $\phi$ defined as 
 
$$\psi(.) = 1-{1\over (1+\phi(.)^\gamma)^k}$$
is a ${\mathcal C}^2$ function in  a neighbourhood of $\bar x$, such that
$$(d-\phi)(x)\geq (d-\phi)(\bar x)=0.$$  This implies that $J^{2,-} d(\bar x
)\neq
\emptyset.$ According to some of the properties of $d$ recalled in the
introduction, on such a point $d$ is differentiable and then $\nabla
\phi(\bar x) = \nabla d(\bar x) $ has modulus
$1$. 

One has 
$$\nabla \psi(x) = {k\gamma d(x)^{\gamma-1}\nabla \phi(x)\over (1+d(x)^\gamma)^{k+1}}$$
and  
$$D^2 \psi = {k\gamma d^{\gamma-2}\over (1+d^\gamma)^{k+2}} \left[
(\gamma-1-(k\gamma+1) d^\gamma)\grad \phi \otimes\grad \phi + d(1 +
d^{\gamma}) D^2\phi\right].$$

We need to prove that  one can choose $k$
large enough in order that 
$$F(x, \nabla \psi, D^2\psi)+ b(x).\nabla \psi|\nabla \psi|^\alpha+c(x)
\psi^{\alpha+1}
\leq
-1.$$

We use Remark \ref{dist}  on the distance function and the following inequalities on symmetric matrices 
$$\mbox{if}\  Y\geq 0 \quad \mbox{then} \quad (X-Y)^+\leq X^+\quad\mbox{ and}\quad  (X-Y)^-\geq X^--Y.$$
Using these with $X = D^2\phi$ and $Y = D\phi\otimes D\phi$, and condition (H2)
we obtain
\begin{eqnarray*}
&F(x,\grad \psi, D^2\psi) 
\leq & \\
& {k^{\alpha+1} \gamma^{1+\alpha}d^{\gamma(\alpha+1)-2-\alpha}\over
(1+d^\gamma)^{k(\alpha+1)+\alpha+2}}
\left[a(\gamma-1 -d^\gamma(k\gamma+1))+(A+a)C_1Nd(1+d^\gamma)\right] & 
\end{eqnarray*}
The function 
${d^{\gamma(\alpha+1)-2-\alpha}\over
(1+d^\gamma)^{k(\alpha+1)+\alpha+2}}$ is decreasing hence it is greater than 
$${K^{\gamma(\alpha+1)-2-\alpha}\over
(1+K^\gamma)^{k(\alpha+1)+\alpha+2}}= C_4.$$ We shall use this later. 

Now we write 
\begin{eqnarray*}
b(x).\nabla \psi|\nabla \psi|^{\alpha} & \leq & |b|_\infty k^{1+\alpha}
\gamma^{1+\alpha} {d^{(\gamma-1)(1+\alpha)}\over
(1+d^\gamma)^{(k+1)(\alpha+1)}}\\
 & \leq &|b|_\infty  {k^{\alpha +1}
\gamma^{1+\alpha}d^{\gamma(\alpha+1)-2-\alpha}\over (1+d^\gamma)^{(k+1)(\alpha+1)+1}}d(1+d^\gamma)
\end{eqnarray*}

We have obtained that there exists a constant $C =C(A,a,|b|_\infty,N)$ such
that 
$$
\displaystyle F(x,\grad\psi, D^2\psi) + b(x).\nabla \psi|\nabla \psi|^\alpha +
c(x)\psi^{1+\alpha}\leq 
$$
$$
\displaystyle {k^{\alpha+1}
\gamma^{1+\alpha}d^{\gamma(\alpha+1)-2-\alpha}\over
(1+d^\gamma)^{k(\alpha+1)+\alpha+2}}
\left[a(\gamma-1
-d^\gamma(k\gamma+1))
+C
d(1+d^\gamma)
\right].
$$
Clearly since $\gamma<1$ we can choose $k$ large enough in order that 
$$[a(\gamma-1
-d^\gamma(k\gamma+1))+Cd(1+d^\gamma)]<{-1\over C_4 k^{1+\alpha} \gamma^{1+\alpha}}<0.$$ 
Then
$$G(x, \psi, \grad \psi, D^2\psi)\leq k^{\alpha+1}
\gamma^{1+\alpha}C_4({-1\over C_4  (k\gamma)^{1+\alpha}})= -1$$ 
which gives the result. 
This ends the proof.

\begin{remark}\label{MM} Clearly if $u$ is the super solution 
constructed in the previous Proposition then  for any $M>0$ and any
$0\leq c_o\leq \left(\frac{M}{|c|_\infty}\right)^{1\over
1+\alpha}$ the function $u_2(x)=M^{1\over 1+\alpha}u(x)+c_o$ is a super
solution of:
$$\left\{\begin{array}{lc}
F(x,\grad u_2,D^2u_2) +b(x).\nabla u_2 |\nabla
u_2|^{\alpha+1}+c(x)u_2^{1+\alpha}\leq -M &  {\rm in}\
\Omega\\ 
u_2=c_o\ & {\rm on}\ \partial\Omega.
\end{array}\right.
$$
\end{remark}

\bigskip
\noindent We are now in a position to solve step 2:

\begin{theorem}\label{propflambda}
Suppose that $F$ satisfies (H1) and (H2), that
$b$ and $c$ are continuous with  $c\leq 0$. 
\begin{enumerate}
\item If $f$ is continuous, bounded and $f\leq 0$ on
$\overline{\Omega}$,   then there exists
$u$  a nonnegative viscosity solution of  
$$\left\{\begin{array}{lc}
F(x,\grad u,D^2u) +b(x).\nabla u |\nabla
u|^{\alpha+1}+c(x)u^{1+\alpha}= f &  {\rm in}\
\Omega\\ u=0\ & {\rm on}\ \partial\Omega.
\end{array}\right.
$$

\item For any  bounded continuous function $f<-M<0$ for some positive
 constant $M$   and any $0\leq
c_o\leq\left(\frac{M}{|c|_\infty}\right)^{1\over 1+\alpha}$
there exists $u$ a non negative solution of
\begin{equation}\label{co}\left\{\begin{array}{lc}
F(x,\grad u,D^2u) +b(x).\nabla u |\nabla
u|^{\alpha+1}+c(x)u^{1+\alpha}= f &  {\rm in}\
\Omega\\ 
u=c_o\ & {\rm on}\ \partial\Omega.
\end{array}\right.
\end{equation}
\end{enumerate}
\end{theorem}
\begin{remark}

In a forthcoming paper,\cite{BD9} we prove existence's results with
general data. 
\end{remark} 
{\bf Proof.}

Let $u_2$ be the  viscosity super solution given in Proposition \ref{propbeta} (see Remark \ref{MM}),  of 
$$G(x,u_2,\nabla u_2 , D^2 u_2)\leq  -|f|_\infty,$$ $u_2=0$ on
$\partial\Omega$. 

 We use Perron's method, see Ishii's paper \cite{I}. We define 
$${\mathcal M}= \{ u\geq 0,\ 0\leq u\leq u_2, u\ {\rm is \ a \
sub solution}\}.$$ 
Let $v(x)=\sup_{u\in{\mathcal M}} u(x)$. We prove that $v$ is
both a sub and a super solution.  

We use the same process as in \cite{BD2} to prove that $v^\star$ is a
sub solution. 
 
We now prove that $v_\star$ is a super solution. If not, there would exist $\bar
x\in \Omega$, $r>0$ and $\varphi\in {\mathcal C}^2(B(\bar x,r)$, with  $\nabla
\varphi(\bar x)\neq 0$, satisfying 
$$0= (v_\star-\varphi)(\bar x) \leq (v_\star-\varphi)(x)$$
on $B(\bar x,r)$,  such that
$$G(\bar x, \varphi(\bar x), \grad \varphi(\bar x),D^2\varphi(\bar x))  >f(\bar x).$$

We prove  that $\varphi(\bar x)< v_2(\bar x)$. If not one would have $\varphi(\bar x) = v_\star(\bar x) = v_2(\bar x)$ 
$$(v_2-\varphi)(x)\geq (v_\star-\varphi)(x)\geq (v_\star-\varphi)(\bar x)
= (v_2-\varphi) (\bar x)=0,$$ hence since $v_2$ is a super solution and
$\varphi$ is a test function for $v_2$ on $\bar x$,
$$G(\bar x, \varphi(\bar x) ,  \grad \varphi(\bar x) ,D^2\varphi(\bar x) )
\leq f(\bar x),$$ a contradiction. 
Then  $\varphi(\bar x)< v_2(\bar x)$. We construct now a sub solution which is 
greater than $v_\star $ and less than $v_2$.

Let $\varepsilon>0$ be such that 
$$G(\bar x , \varphi(\bar x) ,\grad \varphi(\bar x) ,D^2\varphi(\bar
x))
\geq
 f(\bar x)+\varepsilon,$$
and let $\delta$ be such that for $|x-\bar x|\leq \delta$:

$$
\displaystyle|G(\bar x ,\varphi(\bar x) ,\grad \varphi(\bar x)
,D^2\varphi(\bar x) )-G( x , \varphi( x) , \grad
\varphi(x),D^2\varphi(x))| 
$$
$$
\displaystyle+ |f(x)-f(\bar x)|
\leq {\varepsilon\over
4}.
$$
Then 
$$G( x , \varphi( x) , \grad
\varphi(x),D^2\varphi(x)) \geq f(x)+{\varepsilon\over 4}.$$
One can assume that 
$$(v_\star-\varphi)(x)\geq |x-\bar x|^4.$$ 
We take $r< \delta^4$ and such that $0<r< \inf_{|x-\bar x|\leq 
\delta}(v_2(x)-\varphi(x))$, and define 
$$w= \sup (\varphi(x)+r, v_\star)$$
$w$ is LSC as it is the supremum of two LSC functions.
  
One has $w(\bar x) = \varphi(\bar x)+ r$, 
and $w= v_\star$ for $r<|x-\bar x|< \delta$.
 
$w$ is a sub solution, since when $w= \varphi+r$ one can use $\varphi+r$ as a test function,  
and using the continuity of $c$, 
$$G( x , \varphi( x) , \grad
\varphi(x),D^2\varphi(x)) \geq
f+{\varepsilon\over 4}.$$
 Elsewhere $w=v_\star$, hence it is a sub solution.
Moreover $w\geq v_\star$, $w\neq v_\star$ and $w\leq g$. This contradicts
the fact that $v_\star$ is the supremum of the sub solutions.  Using
H\"older regularity we get that $v_\star$ is H\"older.

For the proof of the second statement, it is enough  to remark that $u=c_o$ is a sub solution of 
(\ref{co})  and then proceed as above with $u_2$ the solution, given in Proposition   \ref{propbeta} (see Remark \ref{MM}), of
$$G(x,u_2,\nabla u_2 , D^2 u_2)\leq  -|f|_\infty,$$ $u_2=c_o$ on
$\partial\Omega$. 
This ends the proof.

We now prove an existence result for $\lambda<\bar\lambda$ i.e. step 3.

\begin{theorem}\label{prop2flambda}
Suppose that $F$ ,$b$ and $c$   satisfy the assumptions in
Theorem
\ref{compprnew} and that
$\lambda <\bar\lambda$, 

\begin{enumerate}
\item Suppose that
$f\leq 0$, continuous and bounded, then there exists a nonnegative solution of 
$$\left\{\begin{array}{cc}
F(x, \grad u, D^2 u)+b(x).\nabla u|\nabla u|^\alpha+ (c(x)+\lambda )u^{1+\alpha} = f&\ {\rm in}\
\Omega\\ u=0&{\rm on} \ \partial\Omega
\end{array}\right.
$$

\item For any bounded and continuous function $f\leq-M<0$ and any $0\leq c_o\leq\left(
\frac{M}{|c|_\infty}\right)^{1\over 1+\alpha}$
there exists $u$ a non negative solution of
\begin{equation}\label{co}\left\{\begin{array}{lc}
F(x,\grad u,D^2u) +b(x).\nabla u |\nabla
u|^{\alpha+1}+(c(x)+\lambda)u^{1+\alpha}= f &  {\rm in}\
\Omega\\ 
u=c_o\ & {\rm on}\ \partial\Omega.
\end{array}\right.
\end{equation}
\end{enumerate}
\end{theorem}
{\bf Proof.}
We define a sequence by induction  with $u_1=0$ and $u_{n+1}$ as the solution of

$$\left\{\begin{array}{cl}
 F(x, \grad u_{n+1}, D^2 u_{n+1})+b(x).\nabla u_{n+1}|\nabla
u_{n+1}|^\alpha+(c(x) -|c|_\infty )u_{n+1}^{1+\alpha}= &\\[2ex]
= f-(\lambda+|c|_\infty ) u_n^{1+\alpha}  & {\rm in} \  \Omega,
\\[2ex]
u_{n+1}=0 &  {\rm on} \ \partial \Omega.\ 
\end{array}
\right.
$$
which exists by the previous theorem.
 
The sequence is positive and $u_n$
is increasing, indeed  we can use the comparison
 Theorem  \ref{compprnew}  with  the right hand side
equal to 
$f-(\lambda+|c|_\infty) u_n^{1+\alpha} <0$ and the function 
$c(\phi) = (-c+|c|_\infty) \phi^{1+\alpha}$,
which is nonnegative and increasing with respect to $\phi$. 
We need to prove that the sequence is bounded :  suppose that it is not,
then dividing by
$|u_{n+1}|_\infty^{1+\alpha}$ and defining $w_n = {u_n\over |u_n|_\infty}$, one gets that
$w_n$ satisfies 

\begin{eqnarray*}
 F(x, \nabla w_{n+1}, D^2 w_{n+1}) &+ & b(x).\nabla w_{n+1}|\nabla
w_{n+1}|^\alpha + ( c(x)+\lambda-|c+\lambda|_\infty )w_{n+1}^{1+\alpha} =\\[2ex]
&=&{f\over
|u_{n+1}|_\infty^{1+\alpha}}-|\lambda+c|_\infty {u_n^{1+\alpha}\over
|u_{n+1}|_\infty^{1+\alpha}}.
\end{eqnarray*}
Since the sequence is increasing the right hand side is bounded and  it is greater  than
${f\over |u_{n+1}|_\infty^{1+\alpha}}-(|\lambda+c|_\infty)
w_{n+1}^{1+\alpha}$.  
Then $w_n$ converges to $w$ while ${u_n\over
|u_{n+1}|_\infty}$ converges to $k w$ for  $k= \lim
{|u_n|_\infty^{1+\alpha}\over |u_{n+1}|_\infty^{1+\alpha}}\leq 1$.

One gets, by passing to the limit and using the
compactness result, that the limit function $w$ satisfies
$$F(x, w, \grad w, D^2w)+b(x).\nabla w|\nabla w|^\alpha+ (\lambda+c)
w^{1+\alpha}\geq (1-k)|c+\lambda|_\infty w^{1+\alpha}\geq 0$$ with $w\geq 0$, $|w|_\infty=1$ and $w=0$ on the boundary. This contradicts  the maximum principle (Theorem \ref{maxp}). 
  
We have obtained that the sequence $u_n$ is bounded. Letting $n$ go to
infinity, and using the compactness result (Corollary \ref{comp}), the sequence being in
addition monotone, it converges in its whole  to $u$ which is a 
solution.

 The solution is unique if $f\leq -m<0$ on $\overline{\Omega}$. 
Indeed suppose that $u$ and $v$ are two solutions then 
$v(1+\epsilon)$ is a solution with $f(1+\epsilon)^{1+\alpha}$ in the right
hand side. Since  it is strictly less than $f$ one gets  by the
comparison principle \ref{complambda} that
$v(1+\epsilon)\geq u$
and since $\epsilon$ is arbitrary $v\geq u$. One can of course exchange
$u$ and $v$ and obtain that $u = v$
This ends the proof.

\bigskip

\subsection{Existence result for $\lambda = \bar\lambda$}

We have reached
the final step:

\begin{theorem}\label{eigenfunction}
Let $F$, $b$ and $c$ as in Theorem  \ref{compprnew}. Then, there exists
$\phi>0$ in
$\Omega$ such that $\phi$ is a viscosity solution of 

$$
\left\{\begin{array}{lc}
F(x,\grad \phi,D^2\phi)+b(x).\nabla \phi|\nabla
\phi|^\alpha + (c(x)+ \bar\lambda)
\phi^{1+\alpha} = 0 & {\rm in}\ 
\Omega\\
\phi=0  & {\rm on}\  \partial \Omega.
\end{array}
\right.
$$
Moreover $\phi$ is $\gamma$-H\"older continuous for all $\gamma\in ]0,1[$.
\end{theorem}
{\bf Proof.}
Let $\lambda_n$ be an increasing sequence which converges to $\bar\lambda$. Let $u_n$ be a nonnegative viscosity solution of 
$$
\left\{\begin{array}{lc}
F(x, \grad u_n,D^2u_n)+b(x).\nabla u_n|\nabla  u_n|^\alpha+
(c(x)
+\lambda_n) u_n^{1+\alpha} = -1 & {\rm in}\ 
\Omega\\ u_n=0  & {\rm on}\  \partial \Omega.
\end{array}
\right.$$

By  Theorem \ref{propflambda} the sequence $u_n$ is well defined.
We shall prove that $(u_n)$ is not bounded. Indeed suppose by
contradiction that it is.  Then by the H\"older's estimate and the compactness result (Corollary \ref{comp}), one would have that
a subsequence, still denoted $u_n$, tends uniformly to a nonnegative 
continuous function $u$ which would be a viscosity solution of 
$$
\left\{\begin{array}{lc}
F(x, \grad u,D^2u)+b(x).\nabla u|\nabla u|^\alpha
+(c(x)+ \bar{\lambda}) u^{1+\alpha}= -1& \ {\rm in}\ \Omega\\
u=0&\  {\rm on}\ \partial\Omega.
\end{array}
\right.$$
But then  $u >0$ in $\Omega$ and bounded and one can choose 
$\varepsilon$ small enough that
$$F(x, \grad u,D^2u) +b(x).\nabla u|\nabla u|^\alpha
+(c(x)+\bar\lambda+\varepsilon)
u^{1+\alpha}\leq -1+\varepsilon u^{1+\alpha}\leq 0$$
and this contradicts the definition of  $\bar\lambda$.

We have obtained that  $|u_n|_\infty\rightarrow +\infty$. Then defining
$w_n = {u_n\over |u_n|_\infty}$ one has 
$$
\left\{\begin{array}{lc}
G(x,w_n, \grad w_n,D^2w_n)+ \lambda_n w_n^{1+\alpha} = {-1\over
|u_n|^{1+\alpha}} & \ {\rm in}\ \Omega\\
w_n=0  &\  {\rm on}\ \partial\Omega.
\end{array}
\right.$$
and then extracting as previously a subsequence which
converges uniformly, one gets that there exists $w$, $|w|_\infty=1$ and 

$$
\left\{\begin{array}{lc}
G(x,w, \grad w,D^2w)+\bar \lambda w^{1+\alpha} = 0& \ {\rm in}\ \Omega\\
w=0 &\  {\rm on}\ \partial\Omega.
\end{array}
\right.$$

The boundary condition is given by the uniform convergence.
Clearly $w$ is H\"older continuous, and if $F$ satisfies  the
assumption (H7), then it is also  locally Lipschitz continuous. 
This ends the proof.

\begin{remark}
We have   obtained that 
$\bar\lambda$ is also the supremum of the set 
$$\{ \lambda, \exists \phi>0 \ {\rm on}\  \overline{\Omega}, G(x,\phi,
\grad
\phi,D^2\phi)+ \lambda \phi^{1+\alpha}\leq 0\ \mbox{in the viscosity sense}\}.$$
Indeed, for $\lambda<\bar \lambda$ there exists $v$ which is zero  on
the boundary, 
 such that 
$$G(x,v,Dv,D^2)+ \lambda v^{1+\alpha} = -1.$$
Then using a continuity argument, one
gets that for $\epsilon$ small enough $w=v+\epsilon$ is a supersolution
of 
$$G(x,w, Dw, D^2w)+\lambda w^{1+\alpha} \leq
{-1/2}.$$
\end{remark}

\end{document}